\documentclass[12pt]{article}
\usepackage{amsmath,amsthm,amsfonts,amssymb}

\newtheorem{theorem}{Theorem}[section]
\newtheorem{lemma}[theorem]{Lemma}

\newtheorem{corollary}[theorem]{Corollary}
\newtheorem{prop}[theorem]{Proposition}
\newtheorem{assumption}[theorem]{Assumption}
\newtheorem{remark}[theorem]{Remark}

\newcommand {\IN}{\mathbb{N}}  
\newcommand {\IR}{\mathbb{R}}  
\newcommand {\ImB}{\mathcal{B}}  
\newcommand {\cF}{\mathcal{F}}  
\newcommand {\ImF}{\mathcal{F}}  
\newcommand {\cH}{\mathcal{H}}  
\newcommand {\IcC}{\mathbb{C}}
\newcommand {\ImD}{\mathcal{D}}  

\numberwithin{equation}{section}

\begin{document}

\begin{center}{\bf\Large Polarity of points 
 for Gaussian random fields}
\vskip 16pt
{\bf Robert C.~Dalang}\footnote[1]{Institut de math\'ematiques,
Ecole Polytechnique F\'ed\'erale de Lausanne, Station 8,
CH-1015 Lausanne, Switzerland.
Email: robert.dalang@epfl.ch ~~ Partially supported by the Swiss National Foundation for Scientific Research.}, {\bf Carl Mueller}\footnote[2]{Department of Mathematics, University of Rochester, Rochester, NY  14627, U.S.A. http://www.math.rochester.edu/people/faculty/cmlr ~~ Supported in part by an NSF grant.} and {\bf Yimin Xiao}\footnote[3]{Department of Statistics and Probability, Michigan State University, A-413 Wells Hall, East Lansing, MI 48824, U.S.A. Email: xiaoyimi@stt.msu.edu ~~ Partially supported by NSF grants DMS-1309856 and DMS-1307470.}
\vskip 1in

\end{center}
{\small
\noindent{\bf Abstract.} We show that for a wide class of Gaussian random fields, points are polar in the critical dimension.  Examples of such random fields include solutions of systems of linear stochastic partial differential 
equations with deterministic coefficients, such as the stochastic heat equation or wave equation with space-time white noise, or colored noise in spatial dimensions $k \geq 1$. Our approach builds on a delicate covering argument developed by M.~Talagrand (1995, 1998) for the study of fractional Brownian motion, and uses a harmonizable representation of the solutions of these stochastic pde's.
}
\vskip 16pt

\vskip 12pt


\section{Introduction}

Hitting probabilities are one of the most studied features of stochastic processes.  Given a process $X=(X_t)$ with values in $\IR^d$ and a subset $A$ of $\IR^d$, we say that $X$ hits $A$ if 
\begin{equation*}
   P\{X_t\in A\mbox{ for some $t$}\} >0.  
\end{equation*}
The set $A$ is {\em polar} for $X$ if $P\{X_t\in A\mbox{ for some $t$}\} =0$. When $X$ is a Markov process, potential theory gives a necessary and
sufficient condition for a set to be polar: see \cite{bg68} for an extensive discussion.  One first constructs a potential theory associated to $X$, 
after which it follows that $X$ hits $A$ with positive probability if and only if cap$(A)>0$, where cap$(A)$ is the capacity of $A$ with respect to the potential theory associated to $X$.  

For processes other than Markov processes, and even for Gaussian random fields, results on hitting probabilities are much less complete. One exception is the Brownian sheet, which has specific properties such as independence of increments.  Using these properties, Khoshnevisan and Shi \cite{ks99} have given essentially complete answers about hitting probabilities for the sheet, and the recent work of Dalang, Khoshnevisan, Nualart, Wu, and Xiao \cite{dknwx12} and Dalang and Mueller \cite{dm14} has even settled the issue of multiple points of the Brownian sheet in critical dimensions. 

   Other interesting Gaussian random fields are for instance those obtained as solutions of linear systems of stochastic partial differential equations (spde's). Mueller and Tribe \cite{MT} considered systems of $d$ stochastic heat equations
\begin{equation}\label{she1}
   \frac{\partial u}{\partial t}(t,x) = \frac{\partial^2 u}{\partial x^2}(t,x) + \dot W,
\end{equation}
where $t>0$, $x \in \IR$, $\dot W = \dot W(t,x)$ is an $\IR^d$-valued two-parameter white noise, and the function $u(0,\cdot)$ takes values in $\IR^d$ and is suitably specified. This system of spde's is interpreted in integral form in the framework of Walsh \cite{walshsf}. They show (among other things) that points are polar if and only if $d \geq 6$, so that the critical dimension for hitting points is $d=6$ for the random field $u$ and points are polar in this critical dimension. It turns out that the method of \cite{MT} is quite specific and cannot be extended, for instance, even to the case where the system has deterministic but non-constant coefficients.

  Another case in which the issue of polarity in the critical dimension has been resolved concerns systems of reduced stochastic wave equations (in one spatial dimension) studied by Dalang and Nualart in \cite{DN}. In this case, the critical dimension is $d=4$ and points are polar in this dimension (for linear and nonlinear systems of such equations). This situation is again special, because the natural filtration of the process has the commutation property F4 of Cairoli and Walsh \cite{CW}, which makes it possible to use Cairoli's maximal inequality for multiparameter martingales \cite[Chapter 7.2]{khosh}.
	
For linear and nonlinear systems of stochastic heat and wave equations, there has been much progress in recent years for all dimensions except the critical dimension. A typical result is given in \cite{DKN07,DKN09}. In these papers, the authors establish upper and lower bounds on hitting probabilities of the following type:
$$
   c^{-1} \mbox{Cap}_{d-6-\eta}(A) \leq P\{ u(t,x) \in A \mbox{ for some }(t,x) \in [1,2]^2 \} \leq c \cH_{d-6-\eta}(A),
$$
where Cap denotes Bessel-Riesz capacity, $\cH$ denotes Hausdorff measure, and $\eta >0$. This type of upper and lower bound is also available for systems of heat and wave equations in spatial dimensions $k \geq 1$ (see \cite{DKNhi}), for linear systems of stochastic wave equations in spatial dimensions $k \geq 1$ (see \cite{DSb}), and for nonlinear systems of stochastic wave equations in spatial dimensions $k \in \{1,2,3\}$ (see \cite{DSm}). For a wide class of so-called anisotropic Gaussian random fields $v=(v(x),\, x \in \IR^k)$, Bierm\'e, Lacaux and Xiao \cite{BCX} identified the critical dimension and obtained the following result. Let $\alpha_i$ be the H\"older exponent of the random field when the $i$-th coordinate varies and the others are fixed, and set $Q= \alpha_1^{-1} + \cdots + \alpha_k^{-1}$.  Under certain assumptions, they established the following upper and lower bounds on hitting probabilities: Fix $M > 0$ and a compact set $I \subset \IR^k$. Then there is $0<C<\infty$ such that for every compact set $A\subset B(0,M)$ (the open ball in $\IR^d$ centered at $0$ with radius $M$),
$$
   C^{-1} \mbox{ Cap}_{d-Q}(A) \leq P\{\exists x \in I: v(x) \in A \} \leq C \cH_{d-Q}(A).
$$
This result provides lots of information about hitting probabilities when $d \neq Q$ (see also \cite{xiao09}). However, in the case critical where $d=Q$ and $A = \{z_0\}$ is a single point, these two inequalities essentially reduce to $0 \leq P\{\exists x \in I: v(x) = z_0\} \leq 1$, which is uninformative! Some other references on hitting probabilities for linear systems of spde's include \cite{CT,NV,wu}.

   In order to prove that a set is polar, one typically estimates the probability that the random field visits a small ball, and then one uses a covering argument. When the dimension is strictly larger than the critical dimension, rather simple coverings do the job (typically, the covering is obtained via a deterministic partition of the parameter space). For instance, it is rather straightforward to establish that points are polar for standard Brownian in dimensions $d \geq 3$, but the critical dimension $d=2$ is more difficult to handle (see \cite{khoshe}, for instance).
	
	In order to address the issues of exact Hausdorff measure functions and existence of multiple points for a non-Markovian random field such as fractional Brownian motion, Talagrand introduced a new kind of covering argument in the two important papers \cite{T95,T98}. His idea was to consider balls of different (random) sizes that cover a given point in the parameter space. Having noticed that at a typical point, the local (H\"older-type) regularity is better, with high probability, than what one would expect, he chooses ``good balls" that give a sharp cover of the range of the process, allowing the method to succeed even in the critical dimension. His argument relies on properties of Gaussian processes as well as on certain specific properties of fractional Brownian motion. However, it seems that one of his goals was to develop a method that would extend to other situations, since he states, as one reason for studying fractional Brownian motion, that (ordinary) ``Brownian motion suffers from an over abundance of special properties; and that moving away from these forces to find proofs that rely on general principles, and arguably lie at a more fundamental level."
	
	This paper shows that Talagrand's intuition was correct. Indeed, we have isolated sufficient conditions on an anisotropic Gaussian random field $v=(v(x),\, x \in \IR^k)$, as considered in \cite{BCX, xiao09}, under which it is possible to extend Talagrand's argument and establish polarity of points in the critical dimension: see Assumptions \ref{assump1} and \ref{assump2}. These assumptions are satisfied by many multiparameter Gaussian random fields, for which the H\"older exponents in each parameter may be different. The random fields that we consider are typically nowhere-differentiable (see, e.g. Theorem 3.1 in \cite{xiao97} and Theorem 8.1 in \cite{xiao09}), and this assumption states the existence of particular approximations that are Lipschitz continuous but whose Lipschitz constants have a certain asymptotic growth rate. The main assumption \ref{assump1} is discussed in more detail at the beginning of Section \ref{secobj}. This assumption also leads to an upper bound on the canonical metric associated with the Gaussian random field (see Proposition \ref{prop1}).  
	
	 The first technical effort is to establish Proposition \ref{prop4.1}, which extends an analogous result of Talagrand \cite[Proposition 3.4]{T98} and makes precise the idea that for any $x\in \IR^k$, with high probability, there is a (random) neighborhood of $x$ in which the increments $v(y) - v(x)$ are smaller than expected. With this result in hand, and under the assumption that the process has covariances that have better H\"older regularity than its sample paths (see Assumption \ref{assump2}), which is the case in the examples that we are interested in, we extend the method of Talagrand \cite{T98} and establish polarity of points in the critical dimension $Q$ (see Theorem \ref{thm1}). These results are proved in Sections \ref{secprelim}--\ref{sec5}.
	
	The next step is to show that the two main assumptions are satisfied in a wide class of important examples. We begin with the case of linear systems of stochastic heat equations. In Section \ref{sec6}, we consider first the case of constant coefficients, in spatial dimension 1, with space-time white noise as in \eqref{she1}, and recover the result of Mueller and Tribe \cite{MT}: points are polar for this process in dimension $d=6$. However, essentially the same calculations apply to the case of higher spatial dimensions, with spatially homogeneous noise with covariance given by a Riesz kernel with exponent $\beta \in \, ]0,2[$, so we also obtain polarity of points in the critical dimension $d= (4+2k)/(2-\beta)$ for this case (when this fraction is an integer). The verification of Assumption \ref{assump1} relies on a {\em harmonizable representation} of the solution $u(t,x)$ of the stochastic heat equation: see \eqref{e6.1}: this representation is analogous to the spectral representation of stationary processes (see \cite{doob,ito,yaglom}). It also appears in \cite{balan} and is of independent interest.
	
	 As we mentioned above, the method of Mueller and Tribe was not robust enough to extend to systems of heat equations with deterministic but non-constant coefficients. We examine this situation in Section \ref{sec7}, and we obtain, under the assumption that these coefficients have some smoothness properties (expressed in terms of their Fourier transform: see Assumption \ref{assump24}), polarity of points in the critical dimension. This applies in particular to the case of spatial dimension 1 with space-time white noise, and the critical dimension remains $d=6$.
	
	In Section \ref{sec8}, we turn to linear systems of stochastic wave equations with constant coefficients. Here, we consider both the cases of spatial dimension $k=1$ with space-time white noise, and higher spatial dimensions with spatially homogeneous noise with covariance given by a Riesz kernel with exponent $\beta \in \, ]0,2[$. The stochastic wave equation presents additional difficulties because the fundamental solution is irregular (it is not even a function when $k \geq 3$). This means that Walsh formalism does not apply directly and we use the extension of this theory developed by Dalang \cite{D99}. For the spatial dimension $k=1$ with space-time white noise, we show that points are polar in the critical dimension $d=4$, and in higher spatial dimensions, under the assumption $\beta \in [1,2[$, we obtain polarity of points in the critical dimension $d=2(k+1)/(2-\beta)$ (when this fraction is an integer). 
	
	 The method developed by Talagrand and the extensions presented in this paper can also be applied to the issue of multiple points of Gaussian random fields in critical dimensions, and can also be used to study the same type questions for nonlinear systems of spde's. These topics are the subject of research in progress and we expect to present them in future papers. 

\section{Main assumptions and results}\label{secobj}

   Recall that a {\em white noise based on a measure $\nu$} is a set function $A \mapsto W(A)$ defined on $\ImB(\IR^k)$ with values in $L^2(\Omega,\cF,P)$ such that for each $A$, $W(A)$ is a centered normal random variable with variance $\nu(A)$, and when $A \cap B = \emptyset$, then $W(A \cup B) = W(A) + W(B)$ and $W(A)$ and $W(B)$ are independent. If $W(A)$ is a centered normal random vector with values in $\IR^d$ instead of $\IR$ and covariance matrix $\nu(A)\cdot I_d$ (where $I_d$ denotes the $d\times d$ identity matrix), then we say that $A \mapsto W(A)$ is an {\em $\IR^d$-valued white noise.}

   In order to motivate Assumption \ref{assump1} below, recall that a stationary Gaussian process $(v(t),\, t \in \IR)$ admits a spectral representation of the form $v(t) = \int_\IR f(t-s)\, dW_s$, where $f$ is a function and $(W_s)$ is a Brownian motion. For fixed $t\in \IR$, we can define a white noise by setting $v(A,t) = \int_A f(t-s)\, dW_s$. In many cases, when $f(s)$ is smooth and has appropriate decay as $s \to \pm \infty$, it happens that if $\vert t-s \vert \sim 2^{-n/\alpha}$, for some $\alpha >0$, then $v(t) - v(s)$ is well-approximated by $v([2^n, 2^{n+1}[, t) - v([2^n, 2^{n+1}[, s)$. Even though we will not be dealing with stationary processes, but with non-stationary random fields, it is often possible to construct a process that plays the same role as $v(A,t)$. This is the motivation for Assumption \ref{assump1} below, and this assumption will be verified for the solutions to the spde's that we will consider in Sections \ref{sec6}--\ref{sec8}, as we explain just below. 

   Let $v=(v(x),\, x \in \IR^k)$ be a centered continuous $\IR^d$-valued  Gaussian random field with i.i.d.~components. We write $v(x) = (v_1(x),\dots,v_d(x))$.
	
\begin{assumption} Let $I \subset \IR^k$ be a closed box: $I = \prod_{j=1}^k [c_j,d_j]$, where $c_j < d_j$. Let $I^{(\varepsilon)}$ denote an $\varepsilon$-enlargement of $I$, in Euclidean norm. There is a Gaussian random field $(v(A,x),\, A \in \ImB(\IR_+),\, x \in \IR^k)$ and $\varepsilon_0 >0$ such that:

   (a) for all $x \in I^{(\varepsilon_0)}$, $A \mapsto v(A,x)$ is an $\IR^d$-valued white noise, $v(\IR_+,x) = v(x)$, and when $A$ and $B$ are disjoint, $v(A,\cdot)$ and $v(B,\cdot)$ are independent;
	
(b) there are constants $c_0 \in \IR_+$, $a_0 \in \IR_+$ and $\gamma_j >0$, $j=1,\dots,k$, such that for all $a_0\leq a \leq b \leq +\infty$, $x,y \in I^{(\varepsilon_0)}$,
\begin{equation}\label{aa1}
   \Vert v([a,b[,x) - v(x) - v([a,b[,y) + v(y) \Vert_{L^2} \leq c_0\left[\sum_{j=1}^k a^{\gamma_j}\, \vert x_j - y_j\vert + b^{-1} \right]
\end{equation}
and
\begin{equation}\label{aa2}
   \Vert v([0,a_0],x) - v([0,a_0],y) \Vert _{L^2} \leq c_0 \sum_{j=1}^k  \vert x_j - y_j\vert .
\end{equation}
\label{assump1}
\end{assumption}

	In order to see that the above assumption is satisfied by many solutions of spde's, it is necessary in each case to construct the random field $v(A,x)$. Let us consider for example the solution $v(x)$ of the linear one-dimensional heat equation driven by space-time white noise. Then $\IR^k$ will be replaced by $\IR_+ \times \IR$, and the generic variable $x$ above becomes $(t,x)$. We define
$$
   v(A,t,x) =  \int\!\!\!\int_{\max(\vert \tau\vert^{\frac{1}{4}},\, \vert \xi\vert^{\frac{1}{2}})\in A } e^{-i\xi x} \frac{e^{-i\tau t} - e^{-t\xi^2}}{\xi^2 - i\tau} W(d\tau,d\xi),
$$
Then we will see in Section \ref{sec6} that Assumption \ref{assump1} is satisfied (with the exponents $\gamma_1 = 3$, $\gamma_2 = 1$, that is, $\alpha_1 = 1/4$ and $\alpha_2 = 1/2$), as is Assumption \ref{assump2} below. 

   Define $\alpha_j \in\, ]0,1[$ by the relation
$$
   \gamma_j = \alpha_j^{-1} -1,\qquad\mbox{that is, } \alpha_j = (\gamma_j+1)^{-1},
$$
and define a metric
$$
   \Delta(x,y) = \sum_{j=1}^k \vert x_j - y_j \vert^{\alpha_j}.
$$
Consider also the canonical metric associated with $v$:
$$
    d(x,y) = \Vert v(x) - v(y) \Vert_{L^2}.
$$
It turns out that under Assumption \ref{assump1}, the metric $\Delta$ provides an upper bound on the canonical metric.

\begin{prop} Under Assumption \ref{assump1}, for all $x,y \in I^{(\varepsilon_0)}$ with $\Delta(x,y) \leq \min(a_0^{-1},1)$, we have $d(x,y) \leq 4c_0 \Delta(x,y)$.
\label{prop1}
\end{prop}

\proof Fix $x,y \in I^{(\varepsilon_0)}$. Observe that for any $a \geq a_0$,
\begin{align*}
   d(x,y)
	 & \leq \Vert  v(x) - v([a_0,a[,x) - v(y) + v([a_0,a[,y)  \Vert_{L^2}\\
	    & \qquad + \Vert v([a_0,a[,x) - v([a_0,a[,y) \Vert_{L^2}
\end{align*}
and by Assumption \ref{assump1}(a),
\begin{align*}
 \Vert v([a_0,a[,x) - v([a_0,a[,y) \Vert_{L^2} &\leq \Vert v(x) - v([a,\infty[,x) -v(y) + v([a,\infty[,y)\Vert_{L^2}  \\
   & \qquad + \Vert - v([0,a_0[,x) + v([0,a_0[,y) \Vert_{L^2} .
\end{align*}
Applying Assumption \ref{assump1}(b), we see that
\begin{align}\label{eq2a0}
   d(x,y) \leq 
      c_0 \left[\sum_{j=1}^k (a_0^{\alpha_j^{-1} - 1} + a^{\alpha_j^{-1} - 1}) \vert x_j - y_j\vert + a^{-1} + \sum_{j=1}^k \vert x_j - y_j\vert \right].
\end{align}
By hypothesis, $\max_{j=1,\dots,k}\vert x_j - y_j\vert^{\alpha_j} \leq \Delta(x,y) \leq a_0^{-1}$, so we choose $a \geq a_0$ such that $\max_{j=1,\dots,k}\vert x_j - y_j\vert^{\alpha_j} = a^{-1}$. Notice that
\begin{align}\nonumber
  (a_0^{\alpha_j^{-1} - 1} +  a^{\alpha_j^{-1} - 1}) \vert x_j - y_j\vert &= \left[ \left(a_0\, \vert x_j - y_j\vert^{\alpha_j} \right)^{\frac{1-\alpha_j}{\alpha_j}} +
	\left(a\, \vert x_j - y_j\vert^{\alpha_j} \right)^{\frac{1-\alpha_j}{\alpha_j}}\right] 
		 \vert x_j - y_j \vert^{\alpha_j} \\ \nonumber
	 & \leq 2 \left(a\, \vert x_j - y_j\vert^{\alpha_j} \right)^{\frac{1-\alpha_j}{\alpha_j}} \vert x_j - y_j \vert^{\alpha_j} \\
	 & \leq 2 \vert x_j - y_j \vert^{\alpha_j}
	\label{eq2a1}
\end{align}
by the choice of $a$. Now \eqref{eq2a0} and \eqref{eq2a1} imply that
$$
    d(x,y) \leq c_0 \left[2 \sum_{j=1}^k \vert x_j - y_j \vert^{\alpha_j} + \max_{j=1,\dots,k}\vert x_j - y_j\vert^{\alpha_j} + \sum_{j=1}^k \vert x_j - y_j\vert   \right].
$$
For $\Delta(x,y) \leq 1$, since $0< \alpha_j < 1$, we conclude that
$
   d(x,y) \leq 4c_0 \Delta(x,y).
$
\hfill $\Box$
\vskip 16pt

   A first objective is to prove the following analogue for $v$ of Proposition 3.4 of Talagrand \cite{T98}. 
	
\begin{prop} Let
\begin{equation}\label{eQ}
   Q = \sum_{j=1}^k (\gamma_j + 1) = \sum_{j=1}^k  \frac{1}{\alpha_j}.
\end{equation}
Let Assumption \ref{assump1} hold.  Then there are constants $\tilde K < \infty$ and $\rho >0$ with the following property. Given $0< r_0 < \rho$, for all $x_0 \in I$, we have
\begin{align} 
	P\left\{\exists r \in [r_0^2,r_0]\! :\! \sup_{y:\Delta(y,x_0) < r} \vert v(y) - v(x_0)\vert \leq \tilde K \frac{r}{(\log\log\frac{1}{r})^{1/Q}} \right\} 
	 \!\geq \! 1\! -\! \exp\! \left[-\!\left[\log\frac{1}{r_0} \right]^{\frac{1}{2}} \right]\! .
\label{c2}
\end{align}
\label{prop4.1}
\end{prop}


   When $d(y,x_0) \leq 4 c_0 \Delta(y,x_0) \leq 4 c_0 r$, one expects that $v(y) - v(x_0)$ is of order $r$, so Proposition \ref{prop4.1} states that with high probability, there is a $\Delta$-ball of radius $r$ in which the increments $v(y) - v(x_0)$ are smaller than expected. This proposition is proved in Section \ref{secprop41}.
	
   In order to obtain results on polarity of points, we need an additional assumption.
	
\begin{assumption} Let $I \subset \IR^k$ be a closed box and $\varepsilon_0>0$ be as in Assumption \ref{assump1}. 

   (a) There are constants $c >0$ and $\tilde c >0$ such that for all $x,y \in I^{(\varepsilon_0)}$, $d(x,y) \geq c \Delta(x,y)$ and $\Vert v(x) \Vert_{L^2} \geq \tilde c$.
	
	   (b) 
There is $\rho >0$ with the following property. For $x\in I$, there are $x' \in I^{(\varepsilon_0)}$, $\delta_j \in \,]\alpha_j,1]$, $j=1,\dots,k$, and $C >0$ such that for all $y, \bar y \in I^{(\varepsilon_0)}$ with $\Delta(x,y) \leq 2 \rho$ and $\Delta(x,\bar y) \leq 2 \rho$,
$$
   \left\vert E[(v_i(y) - v_i(\bar y)) v_i(x')]\right\vert \leq C \sum_{j=1}^k \vert y_j - \bar y_j\vert^{\delta_j}.
$$
\label{assump2}
\end{assumption}

\begin{remark} (a) Part (a) in Assumption \ref{assump2} is the lower bound on the canonical metric which completes the upper bound in Proposition \ref{prop1}.

   (b) Part (b) in Assumption \ref{assump2} states that covariances are smoother than what one gets from the Cauchy-Schwarz inequality, H\"older continuity and Proposition \ref{prop1}:
\begin{align*}   
\left\vert E[(v_i(y) - v_i(\bar y)) v_i(x')]\right\vert &\leq \Vert v_i(y) - v_i(\bar y) \Vert_{L^2} \Vert v_i(x') \Vert _{L^2}  
	\leq \Vert v_i(x') \Vert _{L^2} \sum_{j=1}^k \vert y_j - \bar y_j\vert^{\alpha_j}.
\end{align*}   
This will be the case in the examples that we will consider.
\label{remassump4}
\end{remark}

   The main results of this section is the following.
	
\begin{theorem} Let Assumptions \ref{assump1} and \ref{assump2} hold for all sufficiently small boxes. Assume that $Q =d$. Then for any closed box $J$ and for all $z \in \IR^Q$,
$$
   P\{\exists x \in J: v(x) = z\} = 0.
$$
\label{thm1}
\end{theorem}

This theorem is proved in Section \ref{sec5}. 
	
\section{Preliminaries}\label{secprelim}

   Following \cite[Section 2]{T95}, we first set up some estimates that are needed.
	
   Recall the number $Q$ defined in \eqref{eQ}. Let $I \subset \IR^k$ be a closed box such that Assumption \ref{assump1} is satisfied. For $x_0 \in I$, the number of balls in metric $d$ of radius $\varepsilon$ needed to cover the set
$$
   S_r(x_0) = \{x\in \IR^k : \Delta(x, x_0) < r \}
$$
is $\leq N_d(S_r,\varepsilon) = c r^{Q} / \varepsilon^Q$ (indeed, $x \in S_r(x_0)$ implies that $\vert x_j - x_{0,j}\vert < r^{\alpha_j^{-1}}$, so the volume of $S_r(x_0)$ with respect to Lebesgue measure is $\leq c r^{Q}$, and by Proposition \ref{prop1}, the volume of a $d$-ball of radius $\varepsilon$ is $\geq \tilde c \varepsilon^{Q}$). 

\begin{lemma} 
Let $D$ be the diameter (in metric $d$) of a subset $S \subset \IR^k$. There is a universal constant $K_0$ such that, for all $u >0$, we have
\begin{align*}
	P\left\{\sup_{x, y \in S} \vert v(x) - v(y)\vert \geq K_0\left( u + \int_0^D \sqrt{\log N_d(S,\varepsilon)}\, d\varepsilon \right) \right\} 
		\leq \exp \left(-\frac{u^2}{D^2} \right).
\end{align*}
(Note. There is a misprint in \cite[Lemma 2.1]{T95}, where $D$ should be $D^2$.)
\label{lem2.1}
\end{lemma}

\proof This is a consequence of inequality (11.4) p.302 in \cite{LT}, which holds for Gaussian processes with $\psi(x) = e^{(x^2)} - 1$. 
\hfill $\Box$
\vskip 12pt

\begin{lemma} 
There is a constant $K>0$ (depending on $c_0$ in Assumption \ref{assump1}) such that for all $u>0$,
$$
   P\left\{\sup_{x, y \in I} \vert v(x) - v(y)\vert  \leq u \right\} \geq \exp\left(-\frac{1}{Ku^Q} \right).
$$
\label{lem2.2}
\end{lemma}


\proof We use the small ball estimate for Gaussian processes (see \cite[(7.13) p.257]{ledoux} or Lemma 2.2 of \cite{T95}):
$$
   P\left\{\sup_{x, y \in I} \vert v(x) - v(y)\vert \leq u \right\} \geq \exp\left(- \frac{\psi(u)}{K} \right),
$$
where $\psi(u) = u^{-Q}$. Indeed, a ball of radius $\varepsilon$ (in the canonical metric $d$) has volume $\geq \tilde c \varepsilon^{Q}$, so the number of balls (in the canonical metric $d$) of radius $\varepsilon$ needed to cover $I$ is $\leq c_I \varepsilon^{-Q}$.
\hfill $\Box$
\vskip 16pt

\begin{lemma} 
Consider $b > a > 1$, $\varepsilon_0 > r > 0$ and set 
$$
   A =  \sum_{j=1}^k a^{\alpha_j^{-1} -1}\, r^{\alpha_j^{-1}} + b^{-1}. 
$$
There are constants $A_0$, $\tilde K$ and $\tilde c$ (depending on $c_0$ in Assumption \ref{assump1}) such that if $A \leq A_0 r^{}$ and
\begin{equation}\label{b1}
   u \geq \tilde K A \log^{1/2}\left(\frac{r^{}}{A} \right) ,
\end{equation}
then
\begin{align*}
	P\left\{\sup_{x \in S_r(x_0)} \vert v(x) - v(x_0) - (v([a,b],x) - v([a,b],x_0))\vert \geq u \right\} 
	\leq \exp\left(-\frac{u^2}{\tilde cA^2} \right).
\end{align*}
\label{lem3.2}
\end{lemma}


\proof Recall that $S_r = S_r(x_0) = \{x \in \IR^k : \Delta(x,x_0) < r \}$, and set
$$
 \tilde d(x,y) = \Vert v(x) - v(y) - (v([a,b[,x) - v([a,b[,y) \Vert_{L^2}.
$$
Then
$$
   \tilde d(x,y) \leq \Vert v(x) - v(y) \Vert_{L^2} +  \Vert v([a,b[,x) - v([a,b[,y) \Vert_{L^2}.
$$
Since
$$
   v(x) - v(y) = (v([a,b[, x) - v([a,b[,y)) + (v(\IR_+ \setminus [a,b[, x) - v(\IR_+ \setminus [a,b[,y)),
$$
and the two terms on the right-hand side are independent by Assumption \ref{assump1}(a), we see that
$$
   \Vert v([a,b[, x) - v([a,b[,y) \Vert_{L^2} \leq \Vert v(x) - v(y) \Vert_{L^2}.
$$
Finally, 
$$
   \tilde d(x,y) \leq 2 \Vert v(x) - v(y) \Vert_{L^2} \leq 8 c_0 \Delta(x,y)
$$
by Proposition \ref{prop1}. Therefore, for small $\varepsilon >0$, the number of $\varepsilon$-balls (in metric $\tilde d$) needed to cover $S_r(x_0)$ is
$$
   N_{\tilde d}(S_r(x_0), \varepsilon) \leq c \frac{r^{Q}}{\varepsilon^Q}.
$$
For $x \in S_r(x_0)$, $\vert x_j - x_{0,j} \vert \leq r^{\alpha_j^{-1}}$, so by Assumption \ref{assump1}(b),  $\tilde d(x,x_0) \leq c_0 A$, and therefore the diameter $D$ of $S_r(x_0)$ satisfies $D \leq 2 c_0 A$. Assuming that we have chosen the constant $A_0$ and that $A \leq A_0 r$, notice that for $D \leq 2 c_0 A \leq 2 c_0 A_0 r$, there is a constant $\tilde K'$ (depending on $c$ and $c_0 A_0$) such that
$$
   \int_0^D \sqrt{\log N_{\tilde d}(S_r(x_0), \varepsilon)}\, d\varepsilon \leq \tilde K' \int_0^D \sqrt{\log\frac{r^{}}{\varepsilon}}\, d\varepsilon.
$$
Recalling the elementary inequality $\int_x^{+ \infty} u^2 e^{-u^2}\, du \leq C x e^{-x^2}$ for  $x$ large, and using the change of variables $\varepsilon = r e^{-u^2}$ ($r$ fixed), we see that there is a universal constant $K$ such that for all $D>0$ and $r>0$ with $D/r$ sufficiently small (which is the case if $A_0$ is chosen sufficiently small),
$$
   \int_0^D \sqrt{\log\frac{r^{}}{\varepsilon}}\, d\varepsilon \leq K D \sqrt{\log\frac{r}{D}},
$$
so for $D/r$ sufficiently small, 
$$
   \int_0^D \sqrt{\log N_{\tilde d}(S_r(x_0), \varepsilon)}\, d\varepsilon \leq \tilde K' K D \sqrt{\log\frac{r^{}}{D}}.
$$
Let $K_0$ be the universal constant in Lemma \ref{lem2.1}. It follows that
$$
   u \geq  K_0 \left(\frac{u}{2K_0} + \int_0^D \sqrt{\log N_{\tilde d}(S_r(x_0), \varepsilon)}\, d\varepsilon \right)
$$
when 
\begin{equation}\label{b2}
   u \geq 2 K_0 \tilde K' K D \sqrt{\log\frac{r^{}}{D}},
\end{equation}
so by Lemma \ref{lem2.1} (applied to the random field $(v(x) - v([a,b[,x))$), when $u$ satisfies \eqref{b2},
\begin{align*}
  & P\left\{\sup_{x \in S_r(x_0)} \vert v(x) - v(x_0) - (v([a,b],x) - v([a,b],x_0))\vert \geq u \right\} \\
	&\qquad \leq \exp\left(-\frac{(u/(2K_0))^2}{D^2} \right) 
	 \leq \exp\left(-\frac{u^2}{\tilde c A^2} \right).
\end{align*}

   In order to explain \eqref{b1}, notice that
$$
	   D \sqrt{\log\frac{r^{}}{D}} = [f(D)]^{1/2}, \qquad \mbox{ where } f(x) = \frac{x^2}{2} \log\frac{r^{2}}{x^2},
$$
and
$$
	   f'(x) = x \log\frac{r^{2}}{x^2} - \frac{x^2}{2}\, \frac{2}{x} = x \left(\log \frac{r^{2}}{x^2} - 1 \right),
$$
so
$$
	   f'(x) > 0 \qquad \mbox{if}\qquad \frac{r^{2}}{x^2} > e, \qquad\mbox{i.e. } x^2 < \frac{r^{2}}{e}.
$$
Since $D \leq 2c_0 A$,
$$
	   K D \sqrt{\log\frac{r}{D}} \leq  \tilde K'' \left[\frac{A^2}{2} \log\left( \frac{r^{2}}{(2c_0A)^2}\right) \right]^{1/2} \leq \tilde K A \log^{\frac{1}{2}}\left(\frac{r}{A}\right)
$$
provided $(2c_0 A)^2 \leq r^{2}/e$, that is, $A \leq  (2 c_0 e)^{-1} r$, which is the case as long as $A_0$ is sufficiently small and $A \leq A_0 r$. In this case, \eqref{b1} implies \eqref{b2}.
	\hfill $\Box$
	\vskip 16pt

\begin{lemma} 
There is a constant $K$ (depending on $c_0$ in Assumption \ref{assump1}) such that if $0< u < r^{}$, then for all $0<a<b$,
\begin{align*}
   P\left\{\sup_{x \in S_r(x_0)} \vert v([a,b[,x) - v([a,b[,x_0) \vert \leq u\right\} 
	  \geq \exp\left(-K\frac{r^{Q}}{u^Q} \right).
\end{align*}
\label{lem3.5}
\end{lemma}


\proof As in the proof of Lemma \ref{lem2.2}, we note that the number of balls of radius $\varepsilon$ (in the canonical metric of $v(a,b,\cdot,\cdot)$) needed to cover $S_r(x_0)$ is $\leq c \varepsilon^{-Q} r^{Q}$. Applying the same small ball estimate as in the proof of Lemma \ref{lem2.2}, we obtain the desired conclusion.
\hfill $\Box$
\vskip 16pt

\section{Proof of Proposition \ref{prop4.1}}\label{secprop41}

   Fix $U > 1$. Set $r_\ell = r_0 U^{-2\ell}$ and $a_\ell = U^{2\ell-1}/r_0$. Consider the largest integer $\ell_0$ such that
\begin{equation}\label{b6}
   \ell_0 \leq \frac{\log(1/r_0)}{2 \log U}.
\end{equation}
Then for $\ell \leq \ell_0$, we have $r_\ell \geq r_0^2$.

   It suffices to show that for some large constant $K_2$,
\begin{align*}
   & P\left\{\exists 1 \leq \ell \leq \ell_0: \sup_{x \in S_{r_\ell}(x_0)} \vert v(x) - v(x_0) \vert \leq K_2 \frac{r_\ell^{}}{(\log\log\frac{1}{r_\ell})^{1/Q}} \right\}\\
	  &\qquad \geq 1 - \exp\left(- \left(\log\frac{1}{r_0} \right)^{1/2} \right).
\end{align*}

   It follows from Lemma \ref{lem3.5} that for $K_2$ large enough so that $K/K_2^Q \leq 1/4$,
\begin{align} \nonumber
   & P\left\{\sup_{x \in S_{r_\ell}(x_0)} \vert v([a_\ell,a_{\ell+1}[,x) - v([a_\ell,a_{\ell+1}[,x_0) \vert \leq K_2\frac{r_\ell^{}}{(\log\log\frac{1}{r_\ell})^{1/Q}} \right\} \\ \nonumber
	  & \qquad \geq \exp\left(-\frac{K}{K_2^Q} \frac{r_\ell^{Q}}{r_\ell^{Q}} \left(\log\log\frac{1}{r_\ell}\right) \right)
		\geq \exp\left(-\frac{1}{4}  \left(\log\log\frac{1}{r_\ell}\right) \right)\\
		& \qquad = \left(\log \frac{1}{r_\ell}\right)^{-1/4}.
\label{c3}
\end{align}
Thus, by independence of the $v([a_\ell,a_{\ell+1}[, \cdot)$, $\ell = 1,\dots,\ell_0$,
\begin{align}\nonumber
 & P\left\{\exists \ell \leq \ell_0: \sup_{x \in S_{r_\ell}(x_0)} \vert v([a_\ell,a_{\ell+1}[,x) - v([a_\ell,a_{\ell+1}[,x_0) \vert \leq K_2\frac{r_\ell^{}}{(\log\log\frac{1}{r_\ell})^{1/Q}} \right\} \\ \nonumber
&\quad=  1 - \prod_{1\leq \ell \leq \ell_0} \Big(1 - P\Big\{\sup_{x \in S_{r_\ell}(x_0)} \vert v([a_\ell,a_{\ell+1}[,x) - v([a_\ell,a_{\ell+1}[,x_0) \vert \\
   & \qquad\qquad\qquad\qquad\qquad\qquad\qquad \leq K_2\frac{r_\ell^{}}{(\log\log\frac{1}{r_\ell})^{1/Q}}\Big\}\Big).
\end{align}
Apply \eqref{c3} to see that this is greater than
\begin{align}
 1 - \prod_{\ell =1}^{\ell_0} \left[1 - \left[\log\frac{1}{r_\ell}\right]^{-\frac{1}{4}} \right]
\geq 1 - \left[1- \left[\log\frac{1}{r_0^2}\right]^{-\frac{1}{4}} \right]^{\ell_0} 
\geq 1 - \exp\left[- \ell_0 \left[\log\frac{1}{r_0^2}\right]^{-\frac{1}{4}}\right].
\label{4.2}
\end{align}

   Set 
$$
   A_\ell = \sum_{j=1}^k a_\ell^{\alpha_j^{-1}-1} r_\ell^{\alpha_j^{-1}} + a_{\ell+1}^{-1}. 
$$
Notice that $r_\ell a_\ell = U^{-1}$ and $r_\ell a_{\ell+1} = U$. Then
$$
  A_\ell r_\ell^{-1} = \sum_{j=1}^k (a_\ell r_\ell)^{\alpha_j^{-1}-1} + (a_{\ell+1} r_\ell)^{-1} = \sum_{j=1}^k U^{-(\alpha_j^{-1}-1)} + U^{-1} \leq (k+1) U^{-\beta},
$$
with $ \beta = \min(1,\min_{j=1,\dots,k} (\alpha_j^{-1}-1)) >0$ since $\alpha_j < 1$, $j=1,\dots,k$.
Therefore, for $U$ large enough, $A_\ell \leq A_0 r_\ell^{}$, and for $u \geq \tilde K r_\ell^{} U^{-\beta} \sqrt{\log U}$, \eqref{b1} is satisfied (with $A$ there replaced by $A_\ell$ and $r$ by $r_\ell$), so by Lemma \ref{lem3.2},
\begin{align*}
   & P\left\{\sup_{x \in S_{r_\ell}(x_0)} \vert v(x) - v(x_0) - v([a_\ell,a_{\ell+1}[,x) + v([a_\ell,a_{\ell+1}[,x_0) \vert \geq u \right\} \\
	&\qquad \leq \exp\left(-\frac{u^2}{\tilde c A_\ell^2} \right) \leq \exp\left(- \frac{u^2}{c r_\ell^{2}} U^{2\beta} \right).
\end{align*}
Proceeding as in \cite[(4.3)]{T95}, we take $u = K_2 r_\ell^{} (\log\log\frac{1}{r_0})^{-1/Q}$, which is possible provided
$$
   K_2 r_\ell^{} \left(\log\log \frac{1}{r_0} \right)^{-1/Q} \geq \tilde K r_\ell^{} U^{-1} \sqrt{\log U},
$$
that is, provided
\begin{equation}\label{b5}
   U (\log U)^{-1/2} \geq \frac{\tilde K}{K_2} \left(\log\log\frac{1}{r_0} \right)^Q,
\end{equation}
which holds if $U$ is large enough, to get
\begin{align} \nonumber
   & P\left\{\sup_{x \in S_{r_\ell}(x_0)} \vert v(x) - v(x_0) - v([a_\ell,a_{\ell+1}[,x) + v([a_\ell,a_{\ell+1}[,x_0) \vert
	\geq K_2 r_\ell^{} (\log\log\frac{1}{r_0})^{-1/Q} \right\} \\
	&\qquad \leq \exp\left(-\frac{U^{2\beta}}{c(\log\log\frac{1}{r_0})^{2/Q}} \right). 
\label{c1}
\end{align}

   Let
\begin{align*}
	F_\ell &= \left\{ \vert v([a_\ell,a_{\ell+1}[,x) - v([a_\ell,a_{\ell+1}[,x_0)\vert \leq \frac{K_2}{2} \frac{r_\ell^{}}{(\log\log\frac{1}{r_\ell})^{1/Q}}\right\},\\
	G_\ell &= \left\{ \vert v(x) - v(x_0) - v([a_\ell,a_{\ell+1}[,x) + v([a_\ell,a_{\ell+1}[,x_0)\vert \geq \frac{K_2}{2} \frac{r_\ell^{}}{(\log\log\frac{1}{r_\ell})^{1/Q}}\right\}.
\end{align*}
Then
\begin{align} \nonumber
 & P\left\{\exists 1 \leq \ell \leq \ell_0: \sup_{x \in S_{r_\ell}(x_0)} \vert v(x) - v(x_0) \vert \leq K_2\frac{r_\ell^{}}{(\log\log\frac{1}{r_\ell})^{1/Q}} \right\}\\ \nonumber
 &\qquad \geq P\left(\cup_{\ell=1}^{\ell_0} (F_\ell \cap G_\ell^c) \right) \geq P\left(\left(\cup_{\ell=1}^{\ell_0} F_\ell\right) \cap \left(\cap_{\ell=1}^{\ell_0} G_\ell^c \right)\right)  \\ 
  &\qquad \geq  P\left( \cup_{\ell=1}^{\ell_0} F_\ell \right) - P\left( \cup_{\ell=1}^{\ell_0} G_\ell\right).
\label{b3}
\end{align}
By \eqref{4.2},
$$
   P\left( \cup_{\ell=1}^{\ell_0} F_\ell \right) \geq 1 - \exp\left(-\ell_0 \left(\log \frac{1}{r_0^2}\right)^{-1/4} \right),
$$
and by \eqref{c1},
$$
   P\left( \cup_{\ell=1}^{\ell_0} G_\ell\right) \leq \ell_0 \exp\left(-\frac{U^{2\beta}}{c\left(\log\log\frac{1}{r_0} \right)^{2/Q}} \right).
$$

Combining with \eqref{b3}, we get
\begin{align*}
  & P\left\{\exists 1 \leq\ell \leq \ell_0: \sup_{x \in S_{r_\ell}(x_0)} \vert v(x) - v(x_0) \vert \leq K_2\frac{r_\ell^{}}{(\log\log\frac{1}{r_\ell})^{1/Q}} \right\}\\
	&\qquad \geq 1 - \exp\left(- \ell_0 \left(\log\frac{1}{r_0^2}\right)^{-1/4}\right) 
	- \ell_0  \exp\left(-\frac{U^{2\beta}}{c(\log\log\frac{1}{r_0})^{2/Q}} \right).
\end{align*}
Therefore, the proof of \eqref{c2} will be complete provided
\begin{equation}\label{b4}
   \exp\left[- \ell_0 \left[\log\frac{1}{r_0^2}\right]^{-\frac{1}{4}}\right] + \ell_0 \exp\left[\frac{-U^{2\beta}}{c(\log\log\frac{1}{r_0})^{2/Q}} \right] \leq \exp\left[- \left[\log\frac{1}{r_0} \right]^{\frac{1}{2}} \right].
\end{equation}

   Recall the condition \eqref{b5}, and that $\ell_0$ is defined in \eqref{b6}. Therefore, if we set
$$
   U = \left(\log\frac{1}{r_0} \right)^{1/(2\beta)},
$$
then 
for $r_0$ small enough, by \eqref{b6},
$$
   \ell_0 \geq \beta \left(\log\frac{1}{r_0}\right)\, \left(\log\log\frac{1}{r_0}\right)^{-1} \gg 1.
$$
Therefore, the left-hand side of \eqref{b4} is bounded above by
\begin{align*}
	\exp\left[-\frac{\left(\log\frac{1}{r_0} \right)^{3/4}}{c\log\log\frac{1}{r_0}} \right] + \left(1+\log \frac{1}{r_0}  \right) \exp\left[ -\frac{\log \frac{1}{r_0}}{c\left(\log\log\frac{1}{r_0} \right)^{2/Q}}\right] 
	\leq \exp\left[- \left(\log\frac{1}{r_0} \right)^{1/2} \right]
\end{align*}
provided $r_0$ is small enough. This completes the proof of Proposition \ref{prop4.1}. 
\hfill $\Box$
\vskip 16pt

\section{Proof of Theorem \ref{thm1}}\label{sec5}

The main effort in establishing Theorem \ref{thm1} will be to prove the next proposition. 

\begin{prop}\label{prop2} 
Assume that $Q=d$. Let $I$ be a sufficiently small box so that Assumptions \ref{assump1} and \ref{assump2} hold. Let $\varepsilon_0 >0$ be as in Assumption \ref{assump1} and let $\rho$ be as in Assumption \ref{assump2}. Fix $x \in I$, and consider the following (random) subset of $\IR^d$:
$$
   M ( \rho, x) = \{ v(y):\, y \in \IR^k \mbox{ and } \Delta( y , x ) \leq \rho \}.
$$
Then for any $z_{0} \in \IR^{Q}$, $P \{z_{0} \in M (\rho, x) \} = 0.$
\end{prop}

\noindent{\em Proof of Theorem \ref{thm1} (assuming Proposition \ref{prop2})}. Let  $J$ be a closed box and
$$
   M = \{v(y): y \in J \}.
$$
Divide $J$ into a finite union of small boxes $I_\ell$ for which Assumptions \ref{assump1} and \ref{assump2} hold. Let $\rho_\ell > 0$ be given by Assumption \ref{assump2} for $I_\ell$. Since $(S_{\rho_\ell}(x),\, x \in I_\ell)$ is an open cover of $I_\ell$, there are $x_{\ell,1},\dots,x_{\ell,n_\ell} \in I_\ell$ such that $I_\ell \subset \cup_{i=1}^{n_\ell} S_{\rho_\ell} (x_{\ell,i})$. It follows that
$$
   M \subset \cup_\ell \cup_{i =1}^{n_\ell}\, M ( \rho_\ell, x_{\ell,i}),
$$
so for any $z_0 \in \IR^{Q}$,
$$
   P \{ z_{0} \in M \} \leq \sum_\ell \sum_{i=1}^{n_\ell} P \{ z_{0} \in M (\rho_\ell, x_{\ell,i}) \} = 0,
$$
by Proposition \ref{prop2}. It follows that $z_0$ is polar for $v$.
\hfill $\Box$
\vskip 12pt


    We now work towards proving Proposition \ref{prop2}. We proceed as in \cite[Section 3]{T98}. Set
\begin{align*}
B_{\rho} (x) & = \{ y \in \IR^k  : \Delta( y, x)  \leq  \rho\}, \\
B'_{\rho} (x) & =  \{ y \in \IR^k  : \Delta( y, x)  \leq   2 \rho\}.
\end{align*}
Let $x' \in \IR^k$ be given by Assumption \ref{assump2}(b).

   Define two $\IR^{d}$-valued random fields
\begin{align*}
    v^{2} (y)  = E( v(y) \mid  v(x')),\qquad
    v^{1} (y)  =  v(y) - v^{2} (y).
\end{align*}

\begin{remark}{\rm 
(a) Because they are Gaussian and orthogonal, the processes $v^1$ and $v^2$ are independent. Further, $v^1$ is independent of the random vector $v(x')$.

   (b) If we only want to prove that {\em almost all points are polar for $v$} (that is, the range of $v$ has Lebesgue measure zero), then we would not need to introduce the process $v^1$. Here, we will prove that the range of $v^1$ has Lebesgue measure zero, and $v^1$ is quite a good approximation of $v$ (so the range of $v$ also has Lebesgue measure zero). Then we will use the independence of $v^1$ and $v(x')$ to deduce that {\em all} point are polar for $v$.
}
\label{lemindep}
\end{remark}

\begin{lemma} 
The random field $v^2 = (v^2(y),\, y \in B'_{\rho} (x))$ has a continuous version, and there is a finite constant $C$ such that for $y \in B'_{\rho} (x)$ and $\bar y \in B'_{\rho} (x)$,
$$
   \vert v^{2} (y) - v^{2} (\bar y) \vert \leq C   \vert v(x')\vert \sum_{j=1}^k \vert y_{j} - \bar y_{j} \vert^{\delta_j}.
$$
\label{lem4}
\end{lemma}

\proof Let
\begin{equation}\label{a1}
   \alpha(y) = \frac{E ( v_j(y) v_j(x'))}{E ( v_j(x')^{2})},
\end{equation}
where the right-hand side does not depend on $j$. Since the components of $v(y)$ are independent, $v_j^2(y)$ is the orthogonal projection of $v_j(y)$ onto $v_j(x')$, therefore, for $j\in \{1,\dots,d\}$,
\begin{equation}\label{(6)}
   v^{2}_j (y) = \alpha(y) v_j(x'),
\end{equation}
and $v^2(y) = (v^2_1(y),\dots,v^2_d(y))$ is the continuous version of $v^2$. With this version, the conclusion follows from Assumption \ref{assump2}(a) and (b).
\hfill $\Box$
\vskip 16pt

\begin{lemma} 
There is a number $K$ (depending on $d$) such that, for $\varepsilon < 1/3$,
\begin{align*}
 & P\left\{\forall y, \bar y \in \bar I,\ d(y, \bar y) \leq \varepsilon \Rightarrow \vert v(y) - v(\bar y) \vert  \leq K \varepsilon \log^{1/2}\frac{1}{\varepsilon}\right\} 
   \geq 1 - \varepsilon.
\end{align*}
\label{lem2.4}
\end{lemma}

\proof By Assumption \ref{assump2}(a), $I^{(\varepsilon)}$ has finite diameter in the metric $d$. According to \cite[Theorem 6.3.3 p.258]{MR}, there is a finite random variable $Z$ such that, a.s., for all $y, \bar y \in \bar I$,
\begin{align*}
   \vert v(y) - v(\bar y) \vert & \leq Z \int_0^{d(y, \bar y)} \left[\left(\log\frac{1}{\lambda(B_d(y,u))}\right)^{1/2} 
	+ \left(\log\frac{1}{\lambda(B_d(\bar y,u))}\right)^{1/2} \right] du,
\end{align*}
where $\lambda$ denotes Lebesgue on $\bar I$ and $B_d(y,u) $ is the ball in metric $d$ centered at $y$ with radius $u$. Since $d(x,y) \geq c \Delta(x,y)$,
$$
	 \vert v(y) - v(\bar y) \vert \leq Z \int_0^{d(y, \bar y)} \left(\log\frac{1}{u^Q} \right) ^{1/2} du.
$$
Using the elementary inequality
$$
   \int_0^x \left(\log\frac{1}{u} \right)^{1/2} du \leq c_0\, x \log^{1/2}\frac{1}{x},
$$
which is valid for $0<x<x_0$ with $x_0 >0$, and the fact that $x \mapsto x\log^{1/2}(1/x)$ is increasing on $]0,1/e[$, we see that $d(y, \bar y) \leq \varepsilon$ implies 
$$
   \Vert v(y - v(\bar y) \Vert \leq c_0 Z \varepsilon \log^{1/2}(1/\varepsilon), 
$$
and this is $\leq K \varepsilon \log^{1/2}(1/\varepsilon)$ on the event $\{Z \leq K/c_0 \}$. Since $Z$ is finite a.s., this event has probability $\geq 1 - \varepsilon$ if $K$ is large enough.
\hfill $\Box$
\vskip 16pt

   For $p \geq 1, $ consider the random set
\begin{align*} 
	R_{p} = \left\{ y \in B'_{\rho} (x) : \exists r \in [2^{-2p}, 2^{-p} [ \textrm{ with} 
\sup_{\bar{y}: \, \Delta(\bar{y}, y) < r}  \vert v(\bar{y} ) - v (y) \vert \leq K_{2} \frac{r}{(\log \log \frac{1}{r})^{\frac{1}{Q}}} \right\},
\end{align*} 
and the event
$$
   \Omega_{p, 1} = \left\{ \lambda (R_{p}) \geq \lambda (B'_{\rho} (x)) \left(1 - \exp \left ( - \frac{\sqrt{p}}{4} \right) \right) \right\}
$$ 
(here, $\lambda$ denotes Lebesgue measure). Notice that $\Omega_{p, 1}$ can be described as the event ``a large portion of $B'_{\rho} (x)$ consists of points at which $v$ is comparatively smooth." Then
\begin{align*}
(\Omega_{p, 1})^{c} & = \left\{ \lambda (R_{p}) < \lambda (B'_{\rho} (x)) \left( 1 - \exp \left( - \frac{\sqrt{p}}{4} \right)\right) \right\}\\
 & =  \left\{ \lambda ( B'_{\rho} (x) \setminus R_{p}) \geq  \lambda ( B'_{\rho} (x)) \exp \left( - \frac{\sqrt{p}}{4} \right) \right\},
\end{align*}
so by Markov's inequality
\begin{equation}\label{(1)}
 P ( (\Omega_{p,1})^{c}) \leq \frac{E (\lambda ( B'_{\rho} (x) \setminus R_{p}))}{\lambda (B'_{\rho} (x))  \exp \left( - \frac{\sqrt{p}}{4} \right)}\, .
\end{equation}
The numerator is equal to 
$$
   E \left[ \int_{B'_{\rho} (x)} 1_{B'_{\rho} (x) \setminus R_{p}} (y) \, dy\right] =  \int_{B'_{\rho} (x)} P \{ y \in B'_{\rho} (x) \setminus R_{p} \}\, dy.
$$
By the definition of $R_p$ and Proposition \ref{prop4.1} (taking the $\log$ in base 2), for $y \in B'_{\rho} (x)$,
$$
   P \{ y \notin R_{p} \} \leq \exp \left( - \left( \log \frac{1}{2^{-p}} \right)^{\frac{1}{2}} \right) = \exp (-\sqrt{p}),
$$
therefore, by \eqref{(1)},
$$
   P ( ( \Omega_{p, 1} )^{c}) \leq \exp \left( - \frac{3}{4} \sqrt{p} \right).
$$
In particular,
\begin{equation}\label{(2)}
\sum^{\infty}_{p=1} \hspace{0.3cm} P ( (\Omega_{p,1})^{c}) < + \infty.
\end{equation}

   Fix $\beta \in\, ]0, \min(\min_{j=1,\dots,k}(\delta_j \alpha_j^{-1} -1),1)[$ (which is possible since $\delta_j > \alpha_j$, $j=1,\dots,k$) and set
$$
   \Omega_{p,2 } = \{ \vert v( x') \vert \leq 2^{\beta p}\}.
$$
Since $v(x')$ is a normal random vector, $\sum_{p \geq 1} P ( (\Omega_{p, 2})^{c}) < + \infty.$ In addition, on the event $\Omega_{p, 2}$, 
the constant of H\"older continuity of $v^{2}$ is not too large. Indeed, by Lemma \ref{lem4}, for $y \in B'_{\rho} (x)$ and 
$\bar y \in B'_{\rho} (x)$, if $\Delta(y, x) \leq r$ and $\Delta(\bar y, x) \leq r$, then on $\Omega_{p,2}$,
\begin{align*}
    \vert v^{2} (y) - v^{2} (\bar y) \vert \leq C   2^{\beta p} \sum_{j=1}^k \vert y_{j} - \bar y_{j} \vert^{\delta_j}
   \leq \tilde C  \, 2^{\beta p} \, \sum_{j=1}^k r^{\delta_j \alpha_j^{-1}}.
\end{align*}

   If $r \leq 2^{-p},$ then
$$
   r^{\delta_j \alpha_j^{-1}} \, 2^{\beta p} = r^{\delta_j \alpha_j^{-1}} (2^{-p})^{- \beta} \leq r^{\delta_j \alpha_j^{-1}} \, r^{-\beta} = r \, r^{\delta_j \alpha_j^{-1} - 1- \beta},
$$
and $\min_{j=1,\dots,k}(\delta_j \alpha_j^{-1} -1 -\beta) > 0$  by definition of $\beta$. Therefore, there is $K_3 > K_2$ such that on 
$$
   \Omega_{p, 3} \stackrel{\mbox{\scriptsize def}}{=} \Omega_{p, 1} \cap \Omega_{p, 2},
$$ 
for each $y \in R_{p}$, there exists $r \in [2^{-2p}, 2^{-p}]$ such that
\begin{equation}\label{rstar}
   \sup_{\bar{y} :\, \Delta(\bar{y}, y) < r}  \vert v^{1} (\bar{y} ) - v^{1} (y) \vert \leq K_{3} \frac{r}{(\log \log \frac{1}{r})^{\frac{1}{Q}}}.
\end{equation}

   Define an ``anisotropic dyadic cube" of order $\ell$ as a box in $\IR^k$ of the form
$$
   \prod_{j=1}^k [ m_{j}  2^{-\ell \alpha_j^{-1}} , (m_{j}+1)  2^{-\ell \alpha_j^{-1}} ],
$$
where $m_{j} \in \IN.$ For $y \in \IR^k,$ let $C_{\ell} (y)$ denote the anisotropic dyadic cube of order $\ell$ that contains $y.$ This cube is called ``good" if
\begin{equation}\label{(4)}
   \sup_{y, \bar y \in C_{\ell} (y) \cap B_{\rho} (x)} \hspace{0.3cm} \vert v^{1} (y) - v^{1}(\bar y) \vert \leq d_{\ell},
\end{equation}
where
$$
   d_{\ell} = \tilde{K}_{3}\  \frac{2^{-\ell}}{(\log \log 2^{\ell})^{\frac{1}{Q}}}
$$ 
and $\tilde{K}_{3} = k\, K_3$. By \eqref{rstar},  when $\Omega_{p, 3}$ occurs, we can find a family $\mathcal{H}_{1, p}$ of non-overlapping good parabolic dyadic cubes  (they may have intersecting boundaries) of order $\ell \in [p, 2p]$ that covers $R_p$. This family only depends on the random field $v^{1}$. 

   Let $\mathcal{H}_{2, p}$ be the family of non-overlapping dyadic cubes of order $2p$ that meet $B_{\rho} (x)$ but are not contained in any cube of $\mathcal{H}_{1, p}$. For $p$ large enough, these cubes are contained in $B'_{\rho} (x),$ hence in $B'_{\rho} (x) \setminus R_{p}.$ Therefore, when $\Omega_{p, 3}$ occurs, their number is at most $N_{p},$ where
$$
   N_{p} \, 2^{-2pQ} \leq \lambda ( B'_{\rho} (x)) \exp \left( - \frac{\sqrt{p}}{4} \right),
$$
so
\begin{equation}\label{(3)}
   N_{p} \leq C \, 2^{ 2pQ} \exp \left( - \frac{\sqrt{p}}{4} \right),
\end{equation}
where $C$ does not depend on $p$.

   Let $\Omega_{p, 4}$ be the event ``the inequality 
\begin{equation}\label{(5)}
\sup_{y, \bar y \in C} \hspace{0.3cm} \vert v (y) - v(\bar y) \vert \leq K_{4} \, 2^{-2 p} \sqrt{p}
\end{equation}
holds  for each dyadic cube $C$ of order $2p$ of $\IR_{+} \times \IR$ that meets $B_{\rho} (x)$." We choose $K_{4}$ large enough so that $\sum_{p \geq 1} P ( ( \Omega_{p, 4})^{c}) < + \infty :$ this is possible by Lemma \ref{lem2.4}.

   Set $\mathcal{H}_{p} = \mathcal{H}_{1, p} \cup \mathcal{H}_{2, p}.$ This family is well-defined for all $p \geq 1$, and it is a non-overlapping cover  of $B_{\rho} (x)$ (because of how dyadic cubes fit together). Set
$$
  \begin{array}{ll}
   r_{A} = 4  d_{\ell} = 4  \tilde{K}_{3} 2^{- \ell}( \log \ell)^{-\frac{1}{Q}}&\mbox{ if } A \in \mathcal{H}_{1, p} \mbox{ and } A \mbox{ is of order } \ell \in [ p, 2p],\\ [6pt]
   r_{A} = K_{4}  2^{-2 p} \sqrt{p} &\mbox{ if } A \in \mathcal{H}_{2,p}.
	\end{array}
$$

   For each $A \in \mathcal{H}_p$, we pick a distinguished point $p_{A}$ in $A$ (say the lower left corner). Let $B_{A}$ be the ball in $\IR^{d}$ centered at $v (p_{A})$ with radius $r_{A}.$

   Define 
$$
   \Omega_{p} = \Omega_{p, 3} \cap \Omega_{p, 4}.
$$

\begin{lemma} 
Recall that $d=Q$. Let
\begin{equation}\label{f1}
   f(x) = x^{d} \log \log \frac{1}{x}.
\end{equation}
For $p$ large enough, if $\Omega_{p,3}$ occurs, then
$$
   \sum_{A \in \mathcal{H}_p} f(r_{A}) \leq K \lambda (B_{\rho}(x)).
$$
\label{lem2.5}
\end{lemma}

\proof For $A \in \mathcal{H}_{1,p},$
$$
   f(r_{A}) \leq K \left( \frac{2^{- \ell}}{( \log \ell )^{\frac{1}{Q}}} \right)^{d} \log \log 2^{\ell} = K \, 2^{-d \ell}\,  \frac{\log \ell}{(\log \ell)^{d/Q}} = K  \, 2^{-Q \ell}
$$
since $d=Q$, which is the volume of a anisotropic dyadic cube of order $\ell.$

   There is a constant $K_5$ such that, for $p$ large enough and for all $A \in \mathcal{H}_{2,p},$
$$
   f( r_{A} ) \leq K_5 ( 2 ^{-2 p} \sqrt{p} )^{Q}\hspace{0.2cm}  \log (2p).
$$
If $\Omega_{p,3}$ occurs, then by \eqref{(3)}, the total contribution of $\sum_{A \in \mathcal{H}_{2,p}} f(r_{A})$ is bounded by 
$$
   K\, 2^{-2pQ} p^{Q/2} \log (2p) \hspace{0.3cm} 2^{2pQ} \exp \left( - \frac{\sqrt{p}}{4} \right) = p^{Q/2} \log (2p) \exp  \left( - \frac{\sqrt{p}}{4} \right).
$$
Therefore, since the cubes in $\mathcal{H}_{1,p}$ are non-overlapping and intersect $B'_{\rho} (t,x)$, if $\Omega_{p,3}$ occurs, then
$$
   \sum_{A \in \mathcal{H}_p} f(r_{A}) \leq K \lambda (B'_{\rho} (x)) + p^{Q/2} \log (2p) \exp \left( - \frac{\sqrt{p}}{4} \right).
$$
Now $\lambda (B'_{\rho} (x)) \leq 2^{Q} \lambda (B_{\rho} (x)) ,$ and this quantity does not depend on $p$, so the lemma is proved.
\hfill $ \Box$
\vskip 16pt

\begin{lemma} 
Let $\cF_{p}$ be the family of balls $(B_{A},\ A \in \cH_{p})$. For $p$ large enough, on $\Omega_{p}$, $\cF_{p}$ covers $M(\rho, x).$
\label{lem7}
\end{lemma}

\proof Consider $z \in M(\rho, x).$ By definition, there is $y \in B_{\rho} (x)$ such that $v (y) = z.$ Since $\mathcal{H}_{p}$ is a cover of $B_{\rho} (x)$, the point $y$ belongs to a certain cube $A$ of $\mathcal{H}_{p}.$ We will show that $z \in B_{A}$. 
   
	Consider first the case $A \in \mathcal{H}_{1,p}.$ Suppose that $A$ is of order $\ell \in [p,2p]$. By \eqref{(4)},
$$
   \vert v^{1} (p_{A}) - v^{1} (y) \vert \leq d_{\ell}.
$$
Thus, since $\ell \geq p$, on $\Omega_{p,3}$, by Lemma \ref{lem4}, letting $\underline \gamma = \min_{j=1,\dots,k}(\delta_j \alpha_j^{-1} -1 - \beta) >0$,
\begin{align*}
   \vert v (p_{A})- v (y) \vert &\leq d_{\ell} + \vert  v^{2} (v_{A})- v^{2} (y) \vert 
     \leq  d_{\ell} +  C\,2^{\beta p} \sum_{j=1}^k (2^{-\ell})^{\delta_j \alpha_j^{-1}}\\
   & \leq  d_{\ell} + C\,  k\,  2^{\beta p} \, 2^{- \ell\underline \gamma}\, 2^{-\ell(1+\beta)} 
	   \leq d_{\ell} + C\,k\, 2^{- \ell\underline \gamma}\, 2^{-\ell}  \\
   & \leq  2 d_{\ell}
\end{align*}
for $p$ large enough, since $\underline \gamma >0$. Since $v (y) = z$ and $r_{A} = 4 d_{\ell},$ this implies that $z \in B_{A}.$

   Now consider the case $A \in \mathcal{H}_{2}.$ Then on $\Omega_{p,4}$, by \eqref{(5)}, 
$$
   \vert v (p_{A}) -z \vert = \vert v(p_{A}) - v (y) \vert \leq K_{4} \, 2^{-2p} \sqrt{p} = r_{A},
$$
so $z \in B_{A}.$
\hfill $\Box$
\vskip 16pt

\begin{corollary} Almost-surely, the set $M(\rho,x)$ has Lebesgue measure zero: $\lambda(M(\rho,x)) = 0$ a.s.
\label{cor1}
\end{corollary}

\proof For $p$ large enough so that $\Omega_{p}$ occurs, by the definition of $f$ in \eqref{f1} and Lemma \ref{lem2.5},
$$
  \sum_{A \in \mathcal{H}_{p}} r^{d}_{A} \leq \frac{1}{\log p} \sum_{A \in \mathcal{H}_{p}}  f(r_{A}) \leq  \frac{K \lambda (B_{\rho}(x))}{\log p} \to 0
$$
as $p \to + \infty.$ Since the family of balls $(B_{A},\ A \in \cH_{p})$ covers $M(\rho,x)$ by Lemma \ref{lem7}, we conclude that $\lambda(M(\rho,x)) = 0$ a.s.
\hfill $\Box$
\vskip 16pt

\noindent{\em Proof of Proposition \ref{prop2}.} Fix $z_0 \in \IR^Q$. Let $\alpha(y)$ be defined as in \eqref{a1}. Notice that for $\rho$ small enough, $1/2 \leq \alpha(y) \leq 3/2$, and $y \mapsto \alpha(y)$ is H\"older continuous by Assumption \ref{assump2}(b). Define
$$
   v_3(y) = \frac{1}{\alpha(y)}(z_0 - v_1(y)).
$$
Clearly, by \eqref{(6)},
\begin{equation}
   v(y) = z_0 \qquad \Longleftrightarrow \qquad v_3(y) = v(x').
\end{equation}
We are going to check that the range of $v_3$ has Lebesgue measure $0$. Assuming this for the moment, let $f_{v(x')}$ be the probability density function of $v(x')$. Then
\begin{align*}
     P\{z_0 \in M(\rho,x)\} &=   P\{\exists y \in B_{\rho}(x): v_3(y) = v(x')\} \\
		   & = \int_{\IR^Q} dz\, f_{v(x')}(z) P\{\exists y \in B_{\rho}(x): v_3(y) = z\},
\end{align*}
where we have used the fact that $v_1$, hence $v_3$, is independent of $v(x')$ (see Remark \ref{lemindep}(a)). Since the range of $v_3$ has Lebesgue measure $0$, the probability on the right-hand side vanishes for a.a.~$z$, hence the integral is $0$ and $P\{z_0 \in M(\rho,x)\} = 0$, as claimed in Proposition \ref{prop2}.

   It remains to prove that the range of $v_3$ has Lebesgue measure $0$. For $A \in \cH_p$ and $y \in A$,
$$
	v_3(y) - v_3(p_A) = \frac{1}{\alpha(y)} (z_0 - v_1(y)) - \frac{1}{\alpha(p_A)} (z_0 - v_1(p_A)).
$$
Recall that $\alpha$ is H\"older continuous and bounded above and below. If $A \in \cH_{1,p}$ and $A$ is of order $\ell$, then for $p$ sufficiently large, the right-hand side is
$$
    \leq c \sum_{j=1}^k \vert y_j - p_{A,j }\vert^{\delta_j} + 2 d_\ell \leq c \sum_{j=1}^k 2^{-\ell \delta_j \alpha_j^{-1}} +2 d_\ell.
$$
Since $\delta_j \alpha_j^{-1} >1$, $j=1,\dots,k$, this is $\leq 3 d_\ell \leq r_A$. If $A \in \cH_{2,p}$, then for $p$ sufficiently large, the right-hand side is
$$
   \leq c \left(\sum_{j=1}^k 2^{-2p \delta_j \alpha_j^{-1}} + 2 K_4 2^{-2p} \sqrt{2p}\right) \leq  \tilde c\, r_A.
$$
This means that for some constant $\tilde c$, $(B(p_A, \tilde c\,r_A),\, A \in \cH_p)$ covers the range of $v_3$. As in the proof of Corollary \ref{cor1}, we conclude that the Lebesgue measure of $\{v_3(y): \, y \in B_1(x)\}$ is zero.
\hfill $\Box$
\vskip 16pt	
	

\section{Polarity of points for systems of linear heat equations with constant coefficients}\label{sec6}

\noindent Fix $k \geq 1$ and suppose $\beta \in \, ]0, k \wedge2[$ or $k=1=\beta$. Let $\ImD(\IR \times \IR^k)$ denote the space of $C^\infty$-functions with compact support and $(\hat W(\varphi), \ \varphi \in \ImD(\IR \times \IR^k))$ be a spatially homogeneous $\IR^d$-valued Gaussian noise that is white in time, with spatial covariance given by the Riesz kernel $\vert x - y \vert^{-\beta}$, unless $k=1=\beta$, in which case $\hat W$ is space-time $\IR^d$-valued Gaussian white noise based on Lebesgue measure. In both cases, $\hat W(\varphi) =(\hat W_1(\varphi),\dots,\hat W_d(\varphi))$, and the components are independent.

   Recall that in the spatially homogeneous case, the covariance of the noise is informally given by
$$
  E(\hat W_\ell(t,x) \hat W_j(s,y)) = \delta (t-s) \, \vert x-y \vert^{-\beta}\, \delta_{\ell,j},
$$
where $\delta(\cdot)$ denotes the Dirac delta function and $\delta_{\ell,j}$ is the Kronecker symbol. More precisely, for any $C^\infty$-test functions $\varphi$ and $\psi$ with compact support,
$$
E(\hat W_\ell(\varphi) \hat W_j(\psi)) = \delta_{\ell,j} \int_{\IR_+} dr \int_{\IR^k} dy \int_{\IR^k} dz\, \varphi(r,y) \, \vert y-z \vert^{-\beta}\, \psi(r,z).
$$
Using elementary properties of the Fourier transform (see (10) in \cite{D99}), this covariance can also be written
\begin{equation}\label{hnoise}
 E(\hat W_\ell(\varphi) \hat W_j(\psi)) = \delta_{\ell,j} \int_{\IR_+} dr \int_{\IR^k} d\xi\, \vert \xi \vert^{\beta - k}\, \cF_x \varphi(r,\cdot)(\xi) \, \overline{\cF_x\psi(r,,\cdot(\xi)} ,
\end{equation}
where $c_{k,\beta}$ is a constant and $\cF_x \varphi(r,\cdot)(\xi)$ denotes the Fourier transform in the $x$-variable:
$$
   \cF_x \varphi(r,\cdot)(\xi) = \int_{\IR^k} e^{-i \xi \cdot x} \varphi(r,x)\, dx.
$$
This type of noise is discussed for instance in \cite[Section 2]{DKNhi}. Space-time white noise in the case $k=1$ corresponds formally to $\beta = 1$ in \eqref{hnoise}, or, equivalently,
$$
  E(\hat W_\ell(\varphi) \hat W_j(\psi)) = \delta_{\ell,j} \int_{\IR_+} dr \int_{\IR^k} dy \, \varphi(r,y) \psi(r,y).
$$

   Let $\hat v = ( \hat v(t, x), \, t \in \IR_{+}, \, x \in \IR )$ be the mild solution of a linear system of $d$ uncoupled heat equations driven by this space-time white noise: 
\begin{equation}\label{heateq}
\left\lbrace
\begin{array}{rcl}
\frac{\partial}{\partial t} \hat v_{j}(t,x) & = & \Delta \hat v_{j}(t, x) + \dot {\hat W}_{j} (t,x),\qquad j= 1, \dots, d,\\ [8pt]
v (0, x) & = & 0, \qquad x \in \IR^k.
\end{array}
\right.
\end{equation}
Here, $\hat v(t, x) = (\hat v_{1} (t,x), \dots , \hat v_{d} (t,x))$ and $\Delta$ is the Laplacian in the spatial variables. The notion of {\em mild solution} is discussed in \cite[Section 2]{DKNhi} (see also \cite[Chapter 6]{SS}).

\begin{theorem}\label{thm1a}
Suppose $(4+2k)/(2-\beta) =d$. Then points are polar for $\hat v$, that is, for all $z \in \IR^{(4+2k)/(2-\beta)}$, 
$$
   P\{\exists(t,x) \in \, ]0, + \infty[\, \times \IR^k: \hat v (t,x) = z\} = 0.
$$
In particular, in the case where $k=1 = \beta$, $\hat W$ is space-time white noise and $d=6$, then points are polar for $\hat v$.
\end{theorem}

Let $W(d\tau,d\xi)$ be a $\IcC^d$-valued space-time white noise, that is, Re$(W)$ and Im$(W)$ are independent space-time white noises based on Lebesgue measure (Re$(W)$ and Im$(W)$ denote respectively the real and imaginary parts of $W$). In particular,
$$
   E(W_\ell(A)\overline{W_j(B)})= 2 \lambda(A\cap B) \delta_{\ell,j}
$$
(here, $W(A) = (W_1(A),\dots,W_d(A))$). 

   We shall show in the next proposition that the process $(v(t,x),\, (t,x) \in \IR_+ \times \IR^k)$ defined by 
\begin{equation}\label{e6.1}
 v(t,x) = \int_{\IR}\int_{\IR^k} e^{-i\xi\cdot x} \frac{e^{-i\tau t} - e^{-t\vert\xi\vert^2}}{\vert\xi\vert^2 - i\tau} \,\vert \xi\vert ^{(\beta - k)/2}\, W(d\tau,d\xi),
\end{equation}
is a solution of the stochastic heat equation. By analogy with the processes considered in \cite{cambanis}, we call formula \eqref{e6.1} a {\em harmonizable representation} of the solution to \eqref{heateq}. This type of representation also appears in \cite[Section 4]{balan}.

\begin{prop} For $\varphi \in L^2(\IR_+ \times \IR^k,\IcC)$, define
$$
   \tilde W_j(\varphi) = \int_{\IR}\int_{\IR^k} W_j(d\tau,d\xi)\, \vert \xi\vert ^{(\beta - k)/2} \,\cF_{s,y}\varphi(\tau,\xi)  ,
$$
where $\cF_{s,y}$ denotes Fourier transform in the variables $(s,y)$.

   (a) For $j=1,\dots,d$, if $k=1=\beta$, then $\tilde W_j$ is a $\IcC$-valued space-time white noise; otherwise, $\tilde W_j$ is spatially homogeneous noise that is white in time with spatial covariance given by $\vert x-y \vert^{-\beta}$.
	
	 (b) $(v(t,x),\, (t,x)\in \IR_+ \times \IR^k)$ defined in \eqref{e6.1} is a $\IcC$-valued solution of 
\begin{equation}\label{e5.2a}
\left\lbrace
\begin{array}{rcl}
\frac{\partial}{\partial t} v_{j}(t,x) & = & \Delta v_{j}(t, x) + \dot {\tilde W}_{j} (t,x),\qquad j= 1, \dots, d,\\ [8pt]
v (0, x) & = & 0,\qquad x \in \IR^k.
\end{array}
\right.
\end{equation}

   (c) $(\mbox{\rm Re}(v(t,x)), \, (t,x)\in \IR_+ \times \IR^k)$ and $(\hat v(t,x), \, (t,x)\in \IR_+ \times \IR^k)$ have the same law.
\label{prop20}
\end{prop}


\proof (a) Consider first the case $k=1=\beta$. Observe that
\begin{align*}
   E(\tilde W_j(\varphi) \overline{\tilde W_j(\psi)}) & = \int_{\IR} \int_{\IR^k} d\tau d\xi\, \cF_{s,y}\varphi(\tau,\xi) \overline{\cF_{s,y} \psi(\tau,\xi)} 
			= \int_{\IR} \int_{\IR^k} ds dy\, \varphi(s,y) \overline{\psi(s,y)},
\end{align*}
where we have used Plancherel's theorem, so $\tilde W_i$ is a space-time white noise.

   Now consider the case $\beta \in \, ]0,k\wedge 2[$. Then
\begin{align*}
   E(\tilde W_j(\varphi) \overline{\tilde W_j(\psi)}) & = \int_{\IR} d\tau \int_{\IR^k}  \frac{d\xi}{ \vert\xi\vert ^{k- \beta} }\, \cF_{s,y}\varphi(\tau,\xi) \overline{\cF_{s,y} \psi(\tau,\xi)}\\
	    &= \int_{\IR} ds \int_{\IR^k}  dy \int_{\IR^k} dz \, \varphi(s,y) \frac{1}{\vert y - z\vert ^{\beta}} \overline{\psi(s,z)},
\end{align*}	
where we have used again formula (10) in \cite{D99}, and property (a) is established.

   (b) Let $G$ be the fundamental solution of the heat equation. Notice that 
\begin{align*}
   &\int_{\IR\times \IR^k} 1_{[0,t]}(s) G(t-s,x-y) \tilde W_j(ds,dy) \\
	&\qquad = \int_{\IR\times \IR^k} W_j(d\tau,d\xi)\, \cF_{s,y}(1_{[0,t]}(\cdot) G(t-\cdot,x - \cdot))(\tau,\xi) \, \vert \xi\vert ^{(\beta - k)/2}.
\end{align*}
Now, $\cF_{s,y}(1_{[0,t]}(\cdot)G(t-\cdot,x - \cdot))(\tau,\xi)$ is equal to
\begin{align*}
      \cF_s(e^{-i\xi \cdot x} 1_{[0,t]}(\cdot) \overline{\cF_y G(t-\cdot,\cdot)(\xi)})(\tau) 
			& = e^{-i\xi \cdot x}  \cF_s(e^{-(t-\cdot)\vert\xi\vert^2}1_{[0,t]}(\cdot) )(\tau)  \\
			& = e^{-i\xi\cdot x - t\vert\xi\vert^2} \cF_s(e^{s\vert\xi\vert^2} 1_{[0,t]}(s))(\tau).
\end{align*}
The Fourier transform in the $s$-variable is easily calculated and one finds that
\begin{align}\nonumber
   &\int_{\IR\times \IR^k} 1_{[0,t]}(s) G(t-s,x-y) \tilde W_j(ds,dy)\\
	&\qquad = \int_{\IR\times \IR^k} W_j(d\tau,d\xi)\, e^{-i\xi\cdot x} \frac{e^{-i\tau t} - e^{-t\vert\xi\vert^2}}{\vert\xi\vert^2 - i\tau}\, \vert \xi\vert ^{(\beta - k)/2} 
	= v_i(t,x).
\label{p21s1}
\end{align}
By \eqref{e6.1}, $v_j(0,x) = 0$, so, following \cite[Definition 6.1]{SS}, we have checked that $v$ is the (mild) solution of \eqref{e5.2a}, and (b) is proved.

   (c) Set $w = \mbox{Re}(v)$. Then by (b), $w(0,x) = 0$, $w$ satisfies $\frac{\partial w_j}{\partial t} - \Delta w_j  = \mbox{Re}(\tilde W_j(t,x))$. If $k=1=\beta$, then $\mbox{Re}(\tilde W_j)$  is a real-valued space-time white noise such that $E[(\mbox{Re}(\tilde W_j))^2] = \lambda(A)$, and otherwise, $\mbox{Re}(\tilde W_j)$ is a spatially homogeneous noise with the appropriate covariance. This proves (c).
\hfill $\Box$
\vskip 16pt

   Let
\begin{equation}\label{p21s2}
    \alpha_1 = \frac{2-\beta}{4}, \qquad \alpha_2 = \frac{2-\beta}{2} = 2\alpha_1
\end{equation}
(these are the H\"older exponents of $t \mapsto \hat v(t,x)$ and $x \mapsto \hat v(t,x)$, respectively, considered as functions with values in $L^2(\Omega, \cF,P)$), and set
\begin{equation*}
 v(A,t,x) = \int\!\!\!\int_{(\tau,\xi):\, \max(\vert \tau\vert^{\alpha_1},\, \vert \xi\vert^{\alpha_2})\in A } e^{-i\xi\cdot x} \frac{e^{-i\tau t} - e^{-t\vert\xi\vert^2}}{\vert\xi\vert^2 - i\tau} \, \vert \xi\vert ^{(\beta - k)/2}\, W(d\tau,d\xi).
\end{equation*}
Clearly, the random field $(v(A,t,x),\, A \in \ImB(\IR_+),\, (t,x) \in \IR_+ \times \IR^k)$ satisfies Assumption \ref{assump1}(a) (with the generic variable $x \in \IR^k$ replaced by $(t,x) \in \IR_+ \times \IR^k$). In the next lemma, we check Assumption \ref{assump1}(b) (with $a_0 = 0$).

\begin{lemma} 
Let
$$
   \gamma_1 = \alpha_1^{-1} -1 = \frac{2+\beta}{2-\beta},\qquad \gamma_2 = \alpha_2^{-1} -1 = \frac{\beta}{2-\beta}.
$$
There is a universal constant $c_0$ such that for all $0\leq a \leq b$ and $(t_0,x_0)\in\IR_+ \times \IR^k$, $(t,x) \in \IR_+ \times \IR^k$,
\begin{align*}
   &\Vert v([a,b[,t,x) - v(t,x) - v([a,b[,t_0,x_0) + v(t_0,x_0) \Vert_{L^2} \\
	&\qquad \leq c_0\left[a^{\gamma_1}\, \vert t-t_0\vert + a^{\gamma_2}\, \sum_{j=1}^k\vert x_j-x_{0,j}\vert + b^{-1}\right].
\end{align*}
\label{lem3.1}
\end{lemma}

\begin{remark} Lemma \ref{lem3.1} states in particular that for $b = \infty$, $(t,x) \mapsto v(t,x) - v([a,\infty[,t,x)$ is Lipschitz continuous. However, the Lipschitz constants in $t$ and $x$ are of different orders of magnitude, which reflects the $(\alpha_1, \alpha_2)$-H\"older exponents of $(t,x) \mapsto v(t,x)$.
\end{remark}

\proof Let
$$
   v_1(a,t,x) = \int\!\!\!\int_{\max(\vert \tau\vert^{\alpha_1},\, \vert \xi\vert^{\alpha_2})< a } e^{-i\xi x} \frac{e^{-i\tau t} - e^{-t\xi^2}}{\xi^2 - i\tau} \, \vert \xi\vert ^{(\beta - k)/2}\, W(d\tau,d\xi),
$$
$$
   v_2(b,t,x) = \int\!\!\!\int_{\max(\vert \tau\vert^{\alpha_1},\, \vert \xi\vert^{\alpha_2})> b } e^{-i\xi x} \frac{e^{-i\tau t} - e^{-t\xi^2}}{\xi^2 - i\tau} \, \vert \xi\vert ^{(\beta - k)/2}\, W(d\tau,d\xi).
$$
Then
\begin{align}\nonumber
 & v([a,b[,t,x) - v(t,x) - v([a,b[,t_0,x_0) + v(t_0,x_0) \\
 & \qquad = v_1(a,t_0,x_0) - v_1(a,t,x) + v_2(b,t_0,x_0) - v_2(b,t,x).
\label{vab}
\end{align}
Set 
\begin{align*}
   f_1(a,t,x,t_0,x_0) &= E[\vert v_1(a,t,x)- v_1(a,t_0,x_0)\vert^2], \\
	 f_2(b,t,x,t_0,x_0) &= E[\vert v_2(b,t,x)- v_2(b,t_0,x_0)\vert^2].
\end{align*}
We shall estimate these two quantities separately. First, set 
$$
D_1(a) = \{(\tau,\xi)\in \IR \times \IR^k:\max(\vert \tau\vert^{\alpha_1},\, \vert \xi\vert^{\alpha_2})< a \}. 
$$
Then
\begin{align}\nonumber
   & f_1(a,t,x,t_0,x_0)\\  \nonumber
	 &\qquad = d \int\!\!\!\int_{D_1(a) } \left\vert e^{-i\xi \cdot x} \frac{e^{-i\tau t} - e^{-t\vert\xi\vert^2}}{\vert\xi\vert^2 - i\tau} - e^{-i\xi \cdot x_0} \frac{e^{-i\tau t_0} - e^{-t_0\vert\xi\vert^2}}{\vert\xi\vert^2 - i\tau} \right\vert^2 \vert\xi\vert^{\beta - k} d\tau d\xi\\ \nonumber
	&\qquad = d \int\!\!\!\int_{D_1(a)} \left\vert  \frac{e^{-i\tau t} - e^{-t\vert\xi\vert^2} - e^{-i\xi \cdot (x_0 - x) -i\tau t_0} + e^{-t_0\vert\xi\vert^2} e^{-i\xi \cdot (x_0 - x)}}{\vert\xi\vert^2 - i\tau} \right\vert^2 \vert\xi\vert^{\beta - k} d\tau d\xi \\
	&\qquad = d \int\!\!\!\int_{D_1(a)} \frac{\varphi_1(t,x, \tau,\xi)^2 + \varphi_2(t,x, \tau,\xi)^2}{\vert\xi\vert^4 + \tau^2} \vert\xi\vert^{\beta - k} d\tau d\xi,
\label{p23s1}
\end{align}
where
\begin{align*}
   \varphi_1(t,x, \tau,\xi) &= \cos(\tau t) - e^{-t\vert\xi\vert^2} - \cos(\xi \cdot (x_0 - x) + \tau t_0) + e^{-t_0 \vert\xi\vert^2} \cos (\xi\cdot(x_0 - x)),\\
	 \varphi_2(t,x, \tau,\xi) &= -\sin(\tau t) + \sin(\xi \cdot (x_0 - x) + \tau t_0) - e^{-t_0 \xi^2} \sin(\xi \cdot (x_0 - x)).
\end{align*}
Observe that $\varphi_1(t_0,x_0, \tau,\xi) = 0 = \varphi_2(t_0,x_0, \tau,\xi)$, and
\begin{align*}
   \frac{\partial \varphi_1}{\partial t} &= -\tau \sin(\tau t) + \vert\xi\vert^2 e^{-t\vert\xi\vert^2},\\
	\frac{\partial \varphi_1}{\partial x_j} &= -\xi_j \sin(\xi \cdot(x_0 - x) + \tau t_0) + \xi_j e^{-t_0 \xi^2} \sin(\xi \cdot (x_0 - x)),\\
   \frac{\partial \varphi_2}{\partial t} &= -\tau \cos(\tau t), \\
	 \frac{\partial \varphi_2}{\partial x_j} &= -\xi_j \cos(\xi \cdot(x_0 - x) + \tau t_0) + \xi_j e^{-t_0 \vert\xi\vert^2} \cos(\xi \cdot(x_0 - x)).
\end{align*}
Therefore, for $\ell=1,2$,
$$
   \left\vert \frac{\partial \varphi_\ell}{\partial t} \right\vert \leq \vert \tau \vert + \vert\xi\vert^2, \qquad
	  \left\vert \frac{\partial \varphi_\ell}{\partial x_j} \right\vert \leq 2\vert \xi\vert,
$$
and the Mean Value Theorem implies that
\begin{align*}
\vert \varphi_\ell(t,x, \tau,\xi) \vert  &\leq (\vert \tau \vert + \vert\xi\vert^2) \vert t - t_0\vert + 2 \vert \xi \vert \vert x - x_0\vert,
\end{align*}
so
\begin{align}\nonumber
   f_1(a,t,x,t_0,x_0) &\leq d \int\!\!\!\int_{D_1(a)} \left[4 (\tau^2 + \vert\xi\vert^4) (t - t_0)^2 + 8\vert\xi\vert^2\, \vert x - x_0\vert^2\right] 
	\, \frac{\vert\xi\vert^{\beta - k}}{\vert\xi\vert^4 + \tau^2}\, d\tau d\xi\\  \nonumber
	   &\leq 4d \cdot (t - t_0)^2 \int\!\!\!\int_{D_1(a)} \vert\xi\vert^{\beta - k} d\tau d\xi \\ 	
		 & \qquad 
		 + 8d\cdot \vert x-x_0\vert^2 \int\!\!\!\int_{D_1(a)} \frac{\vert\xi\vert^{2+\beta - k}}{\vert\xi\vert^4 + \tau^2} d\tau d\xi.
\label{p23s2}
\end{align}

   For the first integral, pass to polar coordinates $r = \vert \xi \vert$ and use the fact that $\alpha_2 = 2 \alpha_1$ to get
$$
   c_k \int\!\!\!\int_{\max(\vert \tau\vert,r^2) < a^{\alpha_1^{-1}}} r^{\beta-k+k-1} d\tau dr = c_k( A_1 + A_2),
$$
where
$$
   A_1 =  \int\!\!\!\int_{r^2 < \vert \tau\vert < a^{\alpha_1^{-1}}} r^{\beta -1} d\tau dr, \qquad
	 A_2 = \int\!\!\!\int_{\vert \tau\vert < r^2 < a^{\alpha_1^{-1}} } r^{\beta -1} d\tau dr.
$$
Clearly,
$$
   A_1 = \int_0^{a^{\alpha_1^{-1}}} d\tau \int_0^{\sqrt{\tau}} dr\, r^{\beta - 1} = c\,a^{2(2+\beta)/(2-\beta)} = c\, a^{2\gamma_1},
$$
and
$$
   A_2 = \int_0^{a^{\alpha_2^{-1}}} dr\, r^{\beta - 1} \int_0^{r^2}  d\tau = c\, a^{2(2+\beta)/(2-\beta)} = c\, a^{2 \gamma_1}.
$$
We conclude that
\begin{equation}\label{p24s1}
   \int\!\!\!\int_{D_1(a)} \vert\xi\vert^{\beta - k} d\tau d\xi \leq \tilde c a^{2\gamma_1}.
\end{equation}

   For the second integral, pass to polar coordinates $r = \vert \xi \vert$:
\begin{align*}
   \int\!\!\!\int_{D_1(a)} \frac{\vert\xi\vert^{2+\beta - k}}{\vert\xi\vert^4 + \tau^2} d\tau d\xi &=
	   c_k \int_{\max(\vert\tau\vert^{\alpha_1},r^{\alpha_2}) \leq a } \frac{r^{2+\beta -k}}{r^4 + \tau^2} r^{k-1} d\tau dr, 
\end{align*}
then set $w = r^2$ to get
\begin{align*}
   \int_{\max(\vert\tau\vert^{\alpha_1},w^{\alpha_1}) \leq a } \frac{w^{(\beta+1)/2}}{w^2 + \tau^2}  \frac{d\tau dw }{2 \sqrt{w}} 
	 &=  \int_{\max(\vert \tau\vert, w) \leq a^{\alpha_1^{-1}}} \frac{w^{\beta/2}}{w^2 + \tau^2}\, d\tau dw \\
	 &\leq \int_{\max(\vert \tau\vert, w) \leq a^{\alpha_1^{-1}}} \vert (w,\tau) \vert^{\frac{\beta}{2}-2} d\tau dw.
\end{align*}
Pass to polar coordinates $\rho = \vert (w,\tau) \vert$ to see that
\begin{equation}\label{p24s2}
  \int\!\!\!\int_{D_1(a)} \frac{\vert\xi\vert^{2+\beta - k}}{\vert\xi\vert^4 + \tau^2} d\tau d\xi \leq c \int_0^{a^{\alpha_1^{-1}}} d\rho\, \rho^{\frac{\beta}{2}-1} = c a^{2\beta/(2-\beta)} = c a^{2 \gamma_2}.
\end{equation}
We conclude that
$$
   f_1(a,t,x,t_0,x_0) \leq c [ a^{2\gamma_1} (t-t_0)^2 + a^{2\gamma_2} \vert x-x_0\vert^2].
$$

   We now examine $f_2$. Set 
$$
	D_2(b) = \{(\tau, \xi): \max(\vert \tau\vert^{\alpha_1},\, \vert \xi\vert^{\alpha_2})> b \}. 
$$
Notice that as in \eqref{p23s1},
\begin{equation}\label{p25s2}
   f_2(b,t,x,t_0,x_0) = d \int\!\!\!\int_{D_2(b)} \frac{\varphi_1(t,x, \tau,\xi)^2 + \varphi_2(t,x, \tau,\xi)^2}{\vert\xi\vert^4 + \tau^2} \vert\xi\vert^{\beta - k} d\tau d\xi.
\end{equation}
Observing that $\vert \varphi_1 \vert \leq 4$ and $\vert \varphi_2 \vert \leq 3$, we see that
$$
   f_2(b,t,x,t_0,x_0) \leq 25 d \int\!\!\!\int_{D_2(b)} \frac{\vert\xi\vert^{\beta - k}}{\vert \xi\vert^4 + \tau^2} d\tau d\xi.
$$
Let
\begin{align*}
   A_1 &= \{(\tau,\xi): \vert \tau \vert^{\alpha_1} \geq \vert \xi\vert^{\alpha_2} \mbox{ and } \vert \tau \vert^{\alpha_1} > b\}, \\
	A_2 &= \{(\tau,\xi): \vert \tau \vert^{\alpha_1} < \vert \xi\vert^{\alpha_2} \mbox{ and } \vert \xi\vert ^{\alpha_2} > b\} ,
\end{align*}
so that $A_1 \cup A_2 = D_2(b)$. Passing to polar coordinates $\rho = \vert\xi\vert$, notice that
\begin{align}\label{p25s1}
   \int\!\!\!\int_{A_1} \frac{\vert\xi\vert^{\beta - k}}{\vert\xi\vert^4 + \tau^2}\, d\tau d\xi &\leq 4 c_k \int_{b^{\alpha_1^{-1}}}^\infty d\tau \int_0^{\sqrt{\tau}} d\rho\, \frac{\rho^{k-1+\beta -k}}{\tau^2} 
	= c\, b^{2 \frac{\beta - 2}{2 - \beta}} = c\,b^{-2},
\end{align}
and
\begin{align}\label{p25s3}
   \int\!\!\!\int_{A_2} \frac{\vert\xi\vert^{\beta - k}}{\vert\xi\vert^4 + \tau^2}\, d\tau d\xi &\leq c_k \int_{b^{\alpha_2^{-1}}}^\infty d\rho\, \rho^{k-1}  \int_0^{\rho^2} d\tau \frac{\rho^{\beta - k}}{\rho^4}
	=  \tilde c\, b^{2 \frac{\beta - 2}{2 - \beta} } = \tilde c\, b^{-2},
\end{align}
therefore, $f_2(b,t,x,t_0,x_0) \leq c\, b^{-2}$. This proves Lemma \ref{lem3.1}.
\hfill $\Box$
\vskip 16pt

   We now check Assumption \ref{assump2} for the process $\hat v$. In the context of this section, in agreement with Lemma \ref{lem3.1} and Assumption \ref{assump1}(b),
$$
   \Delta((t,x),(s,y)) = \vert t-s \vert^{\frac{2-\beta}{4}} + \vert x-y \vert^{\frac{2-\beta}{2}}.
$$
It is well-known (see \cite[Lemma 4.2]{DKN07}) that for any compact box $I \subset \, ]0,\infty[\, \times \IR$, there is $c >0$ such that for all $(t,x) \in I$ and $(s,y) \in I$,
$$
   \Vert \hat v(t,x) - \hat v(s,y) \Vert_{L^2} \geq c \Delta((t,x),(s,y)).
$$
Further, 
\begin{align}\nonumber
   \Vert \hat v(t,x) \Vert_{L^2}^2 &= d \int_0^t ds \int_{\IR^k} dy \int_{\IR^k} dz \, G(s,y)  \frac{1}{\vert y - z\vert^\beta} G(s,z) \\ \nonumber
	& \geq c^2\int_0^t ds \int_{\IR^k} \frac{d\xi}{\vert \xi\vert^{k - \beta}} \vert \ImF_y G(s,\cdot)(\xi) \vert^2
	=  c^2\int_0^t ds \int_{\IR^k} \frac{d\xi}{\vert \xi\vert^{k - \beta}} e^{-s \vert \xi\vert^2}\\
	 &= c_0 t^{\frac{2-\beta}{2}},
\label{p26s2}
\end{align}
so Assumption \ref{assump2}(a) is satisfied for the box $I$. In the next lemma, we check Assumption \ref{assump2}(b).

\begin{lemma}\label{lem3} 
Let $I \subset \, ]0,\infty[\, \times \IR^k$ be a compact box. Fix $(t,x) \in I$. Let $t' = t - 2 (2\rho)^{\alpha_1^{-1}}$ and $x' = x$ (where $\rho$ is small enough so that $t' >0$). There is a number $C_{1}$ (depending possibly on $\rho$, $\beta$, $k$ and $d$) such that for all
$(s_{1}, y_{1}),\, (s_{2}, y_{2}) \in B'_{\rho} (t,x)$ (the open $\Delta$-ball in $\IR_+ \times \IR^k$ of radius $2\rho$ centered at $(t,x)$), and $j \in \{1, \dots, d\},$
$$
    E[ ( \hat v_{j} (s_{1}, y_{1} ) - \hat v_{j} (s_{2}, y_{2})) \hat v_{j} (t', x')] \leq C_{1} ( \vert s_{1} - s_{2} \vert + \vert y_{1} - y_{2} \vert ).
$$
\end{lemma}

\proof
For $(s, y) \in B'_{\rho} (t,x),$ define
$$
   f(s,y) = E ( \hat v_{j} (s, y) \hat v_{j} (t', x')).
$$

   {\em Case} 1:  $k=1=\beta$. In this case,
$$
   f(s,y) = C \int_{0}^{t'} d r \int_{\IR} d \bar{y} \, G ( s-r, y - \bar{y}) G (t'-r, x'- \bar{y})
$$
(notice that the right-hand side does not depend on $j$). Then
$$
   \frac{\partial f}{\partial y} (s,y) = \int_{0}^{t'} d r \int_{\IR} d \bar{y} \ \frac{\partial G}{\partial y} (s-r, y - \bar{y}) G (t'-r, x'- \bar{y}).
$$
Notice that
$$
   \frac{\partial G}{\partial y} (s-r, y - \bar{y}) = \frac{y - \bar{y}}{s-r} \  G (s-r, y - \bar{y}).
$$
Since $(s, y) \in B'_{\rho} (t,x)$, $s \geq t - (2\rho)^{\alpha_1^{-1}}$, and since $t-t' = 2 (2 \rho)^{\alpha_1^{-1}}$, it follows that for $r \leq t'$, $ s-r  \geq (2\rho)^{\alpha_1^{-1}}$. Therefore $\vert \frac{\partial f}{\partial y} \vert$ is bounded over $B'_{\rho} (t, x)$ (with a bound that depends on $\rho$ but does not depend on $(t,x) \in I$).

   Similarly, since
$$
   \frac{\partial G}{\partial s} (s-r, y - \bar{y}) = - \frac{1}{2} \left( \frac{1}{s-r} + \frac{(y - \bar{y})^{2}}{(s-r)^{2}} \right) G (s-r, y - \bar{y}),
$$
we see that $\vert \frac{\partial f}{\partial s} \vert$ is also bounded over $B'_{\rho} (t, x)$. By the Mean Value Theorem, we conclude that
$$
   \vert f(s_{1}, y_{1} ) - f (s_{2}, y_{2}) \vert \leq C ( \vert s_{1} - s_{2} \vert + \vert y_{1} - y_{2} \vert ),
$$
and this proves the lemma in this case.

   {\em Case} 2: $\beta\in \, ]0,2\wedge k[$. In this case,
\begin{align*}
    f(s,y) &= \int_{0}^{t'} d r \int_{\IR^k} d \bar{y} \int_{\IR^k} d \bar{z} \  G ( s-r, y - \bar{y})  \frac{1}{\vert \bar{y} - \bar{z} \vert^\beta} G (t'-r, x'- \bar{z}) ,
\end{align*}
so
\begin{align*}
  \frac{\partial f}{\partial y_j} (s,y) &= \int_{0}^{t'} d r \int_{\IR^k} d \bar{y} \int_{\IR^k} d \bar{z} \  \frac{\partial G}{\partial y_j} (s-r, y - \bar{y})  \frac{1}{\vert \bar{y} - \bar{z} \vert^\beta} G (t'-r, x'- \bar{z})  \\
	&\leq C \int_{0}^{t'} d r \int_{\IR^k} d \bar{y} \int_{\IR^k} d \bar{z} \  \frac{\partial G}{\partial y_j} (s-r, y - \bar{y})  \frac{1}{\vert \bar{y} - \bar{z} \vert^\beta} G (t'-r, x'- \bar{z}).
\end{align*}
One checks, as above, that $ \frac{\partial f}{\partial y_j} (s,y) $ is bounded (with a bound that depends on $\rho$ but does not depend on $(t,x) \in I$), as is $\frac{\partial f}{\partial s}$, so the conclusion follows as in Case 1.
\hfill $\Box$
\vskip 16pt

\noindent{\em Proof of Theorem \ref{thm1a}.} By Lemma \ref{lem3.1} and the sentences that precede this lemma, for any compact box $I \subset \,]0,\infty[\times \IR^k$, Assumption \ref{assump1} is satisfied for $\mbox{Re}(v)$, with exponents $\gamma_1 = \frac{2+\beta}{2-\beta}$ and $\gamma_j = \frac{\beta}{2-\beta}$, $j=2,\dots,k+1$, so that $\alpha_1 = \frac{2-\beta}{4}$ and $\alpha_j = \frac{2-\beta}{2}$, $j=2,\dots,k+1$. By Lemma \ref{lem3} and the comments that precede this lemma, Assumption \ref{assump2} is satisfied by $\hat v$ (with $\delta_j \equiv 1$), hence by $\mbox{Re}(v)$ by Proposition \ref{prop20}(c). Since $Q = \alpha_1^{-1} + k \alpha_2^{-1} = (4+2k)/(2-\beta) = d$, it follows from Theorem \ref{thm1} that for all $z \in \IR^Q$,  
$$
   P\{\exists (t,x) \in I: \hat v(t,x) = z\} = P\{\exists (t,x) \in I: \mbox{Re}(v(t,x)) = z\} = 0.
$$
Since this holds for all compact rectangles $I \subset \, ]0,\infty[\, \times \IR^k$, Theorem \ref{thm1a} is proved.
\hfill $\Box$
\vskip 16pt

\section{Polarity of points for systems of linear heat equations with nonconstant coefficients}\label{sec7}

 For $j=1,\dots, k$, let $\sigma_j:\IR \times \IR^k \to \IR$ be a continuous functions such that, for all $T \in \IR_+$, there are $0<c_T < C_T< \infty$ such that for all $(t,x) \in [0,T] \times \IR^k$, 
\begin{equation}\label{p27s0}
   c_T \leq \sigma_j(t,x) \leq C_T.
\end{equation}
Let $\hat W$ be as in Section \ref{sec6} and let $\hat v = ( \hat v(t, x), \, t \in \IR_{+}, \, x \in \IR^k )$ be the solution of a linear system of $d$ independent heat equations with deterministic coefficients: 
\begin{equation}\label{heateqv}
\left\lbrace
\begin{array}{rcl}
\frac{\partial}{\partial t} \hat v_{j}(t,x) & = & \Delta \hat v_{j}(t, x) + \sigma_j(t,x) \dot {\hat W}_{j} (t,x),\qquad j= 1, \dots, d,\\ [8pt]
v (0, x) & = & 0, \qquad x \in \IR^k.
\end{array}
\right.
\end{equation}

   Set 
$$
   \tilde G_{t,x}(s,y) = 1_{[0,t]}(s) G(t-s,x-y).
$$
As a consequence of \eqref{p27s0},  in either of the cases $\beta \in \,]0,k\wedge2[$ or $k=1=\beta$,
$$
    \int_\IR d\tau \int_{\IR^k} d\xi\, \vert \xi \vert^{\beta - k} \left\vert\cF_{s,y}( \tilde G_{t,x} \sigma_j)(\tau,\xi)\right\vert^2 
	< \infty.
$$
Indeed, in the case $\beta \in \,]0,k\wedge2[$, for instance, the integral is equal to
$$
   \int_\IR ds \int_{\IR^k} dy_1 \int_{\IR^k} dy_2\,  \tilde G_{t,x}(s,y_1 ) \sigma_j(y_1)\, \frac{1}{\vert y_1 - y_2\vert^\beta}\, G_{t,x}(s,y_2) \sigma_j(y_2),
$$
and then \eqref{p27s0} can be used.

   We also make the following technical assumption on $\sigma_j$. This assumption can be checked for specific choices of $\beta$, $k$ and $\sigma_j$, as in Proposition \ref{prop7.6} below, for instance.

\begin{assumption} (a) $\cF_{s,y}\sigma_j$ is a measure $\mu_j$ with finite total variation.

   (b) Similar to \eqref{p24s2}, for large $a$,
$$
 \int\!\!\!\int_{\IR \times \IR^k} \vert\mu_j\vert(dr,dz) \int\!\!\!\int_{D_1(a)} d\tau d\xi\, \vert \xi \vert^{\beta - k} \frac{\vert \xi - z \vert^2 }{\vert \xi - z \vert^4 + \vert \tau - r\vert^2} \leq c a^{2\gamma_2}.
$$

(c) Similar to \eqref{p25s1} and \eqref{p25s3}, for large $b$,
$$
   \int\!\!\!\int_{\IR \times \IR^k} \vert\mu_j\vert(dr,dz) \int\!\!\!\int_{D_2(b)} d\tau d\xi\, \frac{\vert \xi \vert^{\beta - k} }{\vert \xi - z \vert^4 + \vert \tau - r\vert^2} \leq c b^{-2}.
$$
\label{assump24}
\end{assumption}

\begin{theorem} Suppose that $d = (4+2k)/(2 - \beta)$ and Assumption \ref{assump24} is satisfied. Then points are polar for $\hat v$.
\label{thm25}
\end{theorem}

   Recall from the calculations that led to \eqref{p21s1} that
$$
   \cF_{s,y}\tilde G_{t,x}(\tau,\xi) = e^{-i\xi\cdot x} \frac{e^{-i\tau t} - e^{-t\vert\xi\vert^2}}{\vert\xi\vert^2 - i\tau}.
$$
Define $\tilde W$ as in Proposition \ref{prop20}, and set
$$
   v_j(t,x) = \int_\IR \int_{\IR^k} W_j(d\tau,d\xi) \vert \xi \vert^{(\beta - k)/2} (\cF_{s,y}\tilde G_{t,x} \ast \cF_{s,y}\sigma_j)(\tau,\xi).
$$
 
\begin{prop} The random field $v=(v(t,x) = (v_1(t,x),\dots,v_d(t,x)))$ is the solution of the spde \eqref{heateqv} with $\dot{\hat W}$ replaced by $\dot{\tilde W}$.
\label{prop7.3}
\end{prop} 

\proof Observe that by definition of $\tilde W_j$,
\begin{align*}
  &\int\!\!\!\int 1_{[0,t]}(s)G(t-s,x-y) \sigma(s,y) \tilde W_j(ds,dy)\\
	&\qquad = \int\!\!\!\int W_j(d\tau,d\xi) \cF_{s,y}( \tilde G_{t,x} \sigma_j)(\tau,\xi) \vert \xi \vert^{(\beta - k)/2} 
  = v_j(t,x),
\end{align*}
and $v_j(0,x) = 0$. Therefore, $v$ is the mild solution of \eqref{heateqv} (with $\dot{\hat W}$ replaced by $\dot{\tilde W}$). This completes the proof.
\hfill $\Box$
\vskip 16pt

   Define $\alpha_1$ and $\alpha_2$ as in \eqref{p21s2} and let
$$
   v_j(A,t,x) = \int\!\!\!\int_{\max(\vert \tau\vert^{\alpha_1},\vert \xi\vert^{\alpha_2})\in A } W_j(d\tau,d\xi) \vert \xi \vert^{(\beta - k)/2} (\cF_{s,y}\tilde G_{t,x} \ast \cF_{s,y}\sigma_j)(\tau,\xi)
$$
and $(v(A,t,x) = (v_1(A,t,x),\dots,v_d(A,t,x))$.

\begin{prop} Under Assumption \ref{assump24}(a)--(c), the random field $(v(A,t,x))$ satisfies Assumption \ref{assump1} for any compact box $I \subset \, ]0,\infty[\times \IR^k$.
\label{prop27}
\end{prop}

\proof Assumption \ref{assump1}(a) is clearly satisfied, so we check Assumption \ref{assump1}(b).
   Set
$$
	 v_{1,j}(a,t,x) = \int\!\!\!\int_{D_1(a)} W_j(d\tau,d\xi) \vert \xi \vert^{(\beta - k)/2} (\cF_{s,y}\tilde G_{t,x} \ast \cF_{s,y}\sigma_j)(\tau,\xi).
$$
Define
\begin{align*}
   f_{1,j}(a,t,x,t_0,x_0) &:= E\left(( v_{1,j}(a,t,x) - v_{1,j}(a,t_0,x_0) )^2\right)\\
	   & = \int\!\!\!\int_{D_1(a)} d\tau d\xi\, \vert \xi \vert^{\beta - k}\, \vert ((\cF_{s,y}\tilde G_{t,x} - \cF_{s,y}\tilde G_{t_0,x_0}) \ast \cF_{s,y}\sigma_j)(\tau,\xi) \vert^2,
\end{align*}
and notice that by the Cauchy-Schwarz inequality,
\begin{align}\label{p29b2}
  & \vert ((\cF_{s,y}\tilde G_{t,x} - \cF_{s,y}\tilde G_{t_0,x_0}) \ast \cF_{s,y}\sigma_j)(\tau,\xi) \vert^2 \\ \nonumber
	& \qquad =
	\left\vert \int\!\!\!\int_{\IR \times \IR^k} (\cF_{s,y}\tilde G_{t,x}(\tau - r,\xi-z) - \cF_{s,y}\tilde G_{t_0,x_0}(\tau - r,\xi-z)) \mu_j(dr,dz) \right\vert^2 \\ \nonumber
	&\qquad \leq \vert\mu_j\vert(\IR \times \IR^k) \int\!\!\!\int_{\IR \times \IR^k} \vert \cF_{s,y}\tilde G_{t,x}(\tau - r,\xi-z) 
	- \cF_{s,y}\tilde G_{t_0,x_0}(\tau - r,\xi-z) \vert^2\, \vert\mu_j\vert(dr,dz),
	\nonumber
\end{align}
so
\begin{align*}
   f_{1,j}(a,t,x,t_0,x_0) &\leq C \int\!\!\!\int_{\IR \times \IR^k} \vert\mu_j\vert(dr,dz) \int\!\!\!\int_{D_1(a)} d\tau d\xi\, \vert \xi \vert^{\beta - k} \\
	& \qquad \times \vert \cF_{s,y}\tilde G_{t,x}(\tau - r,\xi-z) - \cF_{s,y}\tilde G_{t_0,x_0}(\tau - r,\xi-z) \vert^2 .
\end{align*}
By \eqref{p23s1} and \eqref{p23s2}, the inner integral is equal to
\begin{align*}
  & d \int\!\!\!\int_{D_1(a)} d\tau d\xi\, \vert \xi \vert^{\beta - k}\ \frac{\varphi_1(t,x,\tau-r,\xi-z)^2 + \varphi_2(t,x,\tau-r,\xi-z)^2}{\vert \xi - z \vert^4 + \vert \tau - r\vert^2} \\
	&\qquad \leq d \int\!\!\!\int_{D_1(a)} d\tau d\xi\, \vert \xi \vert^{\beta - k}\left[4(t-t_0)^2 + 8\frac{\vert x-x_0\vert^2 \vert \xi - z \vert^2 }{\vert \xi - z \vert^4 + \vert \tau - r\vert^2} \right].
\end{align*}
By \eqref{p24s1}, this is
$$
   \leq c_1 (t-t_0)^2 a^{2\gamma_1} + c_2 \vert x-x_0\vert^2 \int\!\!\!\int_{D_1(a)} d\tau d\xi\, \vert \xi \vert^{\beta - k} \frac{\vert \xi - z \vert^2 }{\vert \xi - z \vert^4 + \vert \tau - r\vert^2}.
$$
This establishes in particular \eqref{aa2} for any $a_0 \geq 0$.

   By Assumption \ref{assump24}(a) and (b), we conclude that for large $a$,
\begin{align} \nonumber
   f_1(a,t,x,t_0,x_0) &:= E(\vert v_{1}(a,t,x) - v_{1}(a,t_0,x_0) \vert^2) \\ 
	& \leq c_1\,a^{2\gamma_1} (t-t_0)^2  + c_2\, a^{2\gamma_2} \vert x-x_0\vert^2 .
\label{p29b}
\end{align}

   Set
$$
	 v_{2,j}(a,t,x) = \int\!\!\!\int_{D_2(b)} W_j(d\tau,d\xi) \vert \xi \vert^{(\beta - k)/2} (\cF_{s,y}\tilde G_{t,x} \ast \cF_{s,y}\sigma_j)(\tau,\xi).
$$
Then
\begin{align*}
   f_{2,j}(b,t,x,t_0,x_0) &:= E\left(( v_{2,j}(b,t,x) - v_{2,j}(b,t_0,x_0) )^2\right)\\
	   & = \int\!\!\!\int_{D_2(b)} d\tau d\xi\, \vert \xi \vert^{\beta - k} \vert ((\cF_{s,y}\tilde G_{t,x} - \cF_{s,y}\tilde G_{t_0,x_0}) \ast \cF_{s,y}\sigma_j)(\tau,\xi) \vert^2.
\end{align*}
Using the Cauchy-Schwarz inequality as in \eqref{p29b2}, we find that
\begin{align*}
   f_{2,j}(b,t,x,t_0,x_0) &\leq C \int\!\!\!\int \vert\mu_j\vert(dr,dz) \int\!\!\!\int_{D_2(b)} d\tau d\xi\, \vert \xi \vert^{\beta - k} \\
	& \qquad \times \vert \cF_{s,y}\tilde G_{t,x}(\tau - r,\xi-z) - \cF_{s,y}\tilde G_{t_0,x_0}(\tau - r,\xi-z) \vert^2 
\end{align*}
and by \eqref{p25s2}, the inner integral is equal to 
\begin{align*}
  & \int\!\!\!\int_{D_2(b)} d\tau d\xi\, \vert \xi \vert^{\beta - k}\ \frac{\varphi_1(t,x,\tau-r,\xi-z)^2 + \varphi_2(t,x,\tau-r,\xi-z)^2}{\vert \xi - z \vert^4 + \vert \tau - r\vert^2} \\
	&\qquad \leq 25 \int\!\!\!\int_{D_2(b)} d\tau d\xi\, \frac{\vert \xi \vert^{\beta - k}}{\vert \xi - z \vert^4 + \vert \tau - r\vert^2}.
\end{align*}
By Assumption \ref{assump24}(c), for large $b$,
$$
   f_2(b,t,x,t_0,x_0) := E(\vert v_{2}(b,t,x) - v_{2}(b,t_0,x_0) \vert^2) \leq c b^{-2}.
$$
Putting this together with \eqref{p29b}, we conclude that Assumption \ref{assump1}(b) is satisfied.
\hfill $\Box$
\vskip 16pt

\begin{lemma} $(\hat v(t,x))$ defined in \eqref{heateqv} satisfies Assumption \ref{assump2} for any compact box $I \subset \, ]0,\infty[\times \IR^k$.
\label{lem28}
\end{lemma}

\proof Observe that
\begin{align*}
   \Vert \hat v(t,x) \Vert_{L^2}^2 &= d \int_0^t ds \int_{\IR^k} dy \int_{\IR^k} dz \, G(s,y) \sigma(s,y) \frac{1}{\vert y - z\vert^\beta} G(s,z) \sigma(s,z)\\ 
	 &\geq c_T^2 \int_0^t ds \int_{\IR^k} dy \int_{\IR^k} dz \, G(s,y) \frac{1}{\vert y - z\vert^\beta} G(s,z)
	 \geq c_0 t^{\frac{2 - \beta}{2}}
\end{align*}
by the same calculation as in \eqref{p26s2}, so Assumption \ref{assump2}(a) is satisfied.

   Since $\sigma$ is bounded above by \eqref{p27s0}, the proof of Assumption \ref{assump2}(b) follows the proof of Lemma \ref{lem3}. 
\hfill $\Box$
\vskip 16pt

\noindent{\em Proof of Theorem \ref{thm25}.} By Proposition \ref{prop27}, Assumption \ref{assump1} is satisfied for $\mbox{Re}(v)$, with exponents $\alpha_1 = \frac{2 - \beta }{4}$ and $\alpha_2 = \frac{2 - \beta }{2}$,  for any compact box $I \subset \, ]0,\infty[\times \IR^k$. By Lemma \ref{lem28}, Assumption \ref{assump2} is satisfied for $\hat v$. Since $\mbox{Re}(v)$ and $\hat v$ have the same law by Proposition \ref{prop7.3}, the conclusion follows from Theorem \ref{thm1}.
\hfill $\Box$
\vskip 16pt

\noindent{\em Sufficient conditions for Assumption \ref{assump24}(b) and (c)}
\vskip 12pt

\begin{prop} Suppose that $k=1=\beta$, and $\cF_{s,y}\mu_j$ is a measure with compact support and finite total variation. Then Assumption \ref{assump24} is satisfied.
\label{prop7.6}
\end{prop}

\proof It is clear that Assumption \ref{assump24}(a) holds. We check Assumption \ref{assump24}(b). Note that $k=1=\beta$, so $\alpha_1 = \frac{1}{4}$ and $\alpha_2 = \frac{1}{2}$, and observe that
\begin{align*}
	\int\!\!\!\int_{D_1(a)} d\tau d\xi\, \frac{\vert \xi - z \vert^2 }{\vert \xi - z \vert^4 + \vert \tau - r\vert^2}  
	& \leq \int\!\!\!\int_{\max(\vert \tau - r\vert^{\frac{1}{4}},\vert \xi - z \vert^{\frac{1}{2}}) \leq a + r^{\frac{1}{4}}+z^{\frac{1}{2}}} d\tau d\xi\, \frac{\vert \xi - z \vert^2 }{\vert \xi - z \vert^4 + \vert \tau - r\vert^2} \\ 
	&	= \int\!\!\!\int_{D_1(a+ r^{\frac{1}{4}}+z^{\frac{1}{2}})} d\tau d\xi\, \frac{\vert \xi \vert^2 }{\vert \xi  \vert^4 + \vert \tau \vert^2}.
\end{align*}
By \eqref{p24s2} in the case $k=1=\beta$ (so $\gamma_2 = 1$), we conclude that this integral is $\leq c (a+ r^{\frac{1}{4}}+z^{\frac{1}{2}})^2$, and therefore
\begin{align}\nonumber
  & \int\!\!\!\int_{\IR \times \IR} \vert\mu_i\vert(dr,dz) \int\!\!\!\int_{D_1(a)} d\tau d\xi\,  \frac{\vert \xi - z \vert^2 }{\vert \xi - z \vert^4 + \vert \tau - r\vert^2} \\ \nonumber
	& \qquad \leq c \int\!\!\!\int_{\IR \times \IR} \vert\mu_i\vert(dr,dz) \left(a+ r^{\frac{1}{4}}+z^{\frac{1}{2}}\right)^2 
	= c a^2 \int\!\!\!\int_{\IR \times \IR} \vert\mu_i\vert(dr,dz) \left(1+ \frac{r^{\frac{1}{4}}}{a}+\frac{z^{\frac{1}{2}}}{a}\right)^2 \\
	&\qquad \leq c a^2 \int\!\!\!\int_{\IR \times \IR} \vert\mu_i\vert(dr,dz) \left(1+ r^{\frac{1}{4}}+z^{\frac{1}{2}}\right)^2
\label{eq7.5}
\end{align}
provided $a\geq 1$, and the integral is finite under the assumptions of this proposition. This establishes Assumption \ref{assump24}(b).

   We now check Assumption \ref{assump24}(c) in the case $k=1=\beta$. Use the change of variables $s=\tau - r$, $y = \xi - z$ to see that
\begin{align*}
   \int\!\!\!\int_{D_2(b)} d\tau d\xi\, \frac{1}{\vert \xi - z \vert^4 + \vert \tau - r\vert^2} 
		& \leq \int\!\!\!\int_{\max(\vert s+r\vert^{\frac{1}{4}}, \vert y+z\vert^{\frac{1}{2}}) > b} dsdy\, \frac{1}{\vert y \vert^4 + \vert s \vert^2} \\
		& \leq \int\!\!\!\int_{D_2(\psi(b,r,z))} d\tau d\xi\, \frac{1}{\vert \xi \vert^4 + \vert \tau \vert^2},
\end{align*}
where $\psi(b,r,z) = \min(\left\vert b^4 - \vert r \vert\right\vert^{1/4}, \left\vert b^2 - \vert z\vert\right\vert^{1/2})$.

   By \eqref{p25s1} and \eqref{p25s3}, we conclude that
\begin{align} 
	\int\!\!\!\int_{\IR \times \IR} \vert\mu_i\vert(dr,dz) \int\!\!\!\int_{D_2(b)} d\tau d\xi\, \frac{1}{\vert \xi - z \vert^4 + \vert \tau - r\vert^2} 
	\leq c \int\!\!\!\int_{\IR \times \IR} \vert\mu_i\vert(dr,dz)  (\psi(b,r,z))^{-2},
\label{p33.7.5}
\end{align}
and clearly,
\begin{equation*}
   \psi(b,r,z) = b \min\left(\left\vert1 - \frac{\vert r \vert}{b^4}\right\vert^{1/4}, \left\vert 1 - \frac{\vert z\vert}{b^2}\right\vert^{1/2}\right).
\end{equation*}
If $b$ is large enough so that the inequalities $\vert z \vert \leq b^2/2$ and $\vert r\vert \leq b^4/2$ are satisfied for all $(r,z)$ in the support of $\mu_i$, then $\psi(b,r,z) \geq b/2$, and so the right-hand side of \eqref{p33.7.5} is $\leq 4c\, b^{-2}$. This establishes Assumption \ref{assump24}(c).
\hfill $\Box$
\vskip 16pt

\begin{corollary} Suppose that $d=6$, $k=1$, $\hat W$ is space-time white noise and $\cF_{s,y}\mu_i$ is a measure with compact support and finite total variation. Then points are polar for the solution $(\hat v(t,x))$ of the stochastic heat equation \eqref{heateqv} with nonconstant deterministic coefficients $\sigma_i$.
\end{corollary}

\proof This is an immediate consequence of Theorem \ref{thm25} (with $\beta = 1$) and Proposition \ref{prop7.6}.
\hfill $\Box$
\vskip 16pt


\section{Polarity of points for systems of linear wave equations with constant coefficients}\label{sec8}

Fix $k\geq 1$ and $\beta \in \,]0,k\wedge 2[$ or $k=1 = \beta$, and let $\hat W$ be spatially homogeneous $\IR^d$-valued Gaussian noise as in the beginning of Section \ref{sec6}. We assume that
\begin{equation}\label{ebeta}
   \beta \geq 1.
\end{equation}
Let $\hat v$ be the solution of the stochastic wave equation in spatial dimension $k$ driven by $\hat W$:
$$
\left\lbrace
\begin{array}{rcl}
\frac{\partial^2}{\partial t^2} \hat v_{j}(t,x) & = & \Delta \hat v_{j}(t, x) + \dot {\hat W}_{j} (t,x),\qquad j= 1, \dots, d,\\ [8pt]
\hat v (0, x) & = & 0, \qquad \frac{\partial}{\partial t} \hat v(0,x) = 0, \qquad x \in \IR^k.
\end{array}
\right.
$$

\begin{theorem} Suppose $k=1=\beta$ or $1 < \beta < k \wedge 2$, and $d = \frac{2(k+1)}{2-\beta}$. Then points are polar for $\hat v$, that is, for all $z \in \IR^{d}$, 
$$
   P\{\exists(t,x) \in \, ]0, + \infty[\, \times \IR^k: \hat v (t,x) = z\} = 0.
$$
In particular, in the case where $k=1 = \beta$, $\hat W$ is space-time white noise and $d=4$, then points are polar for $\hat v$.
\label{thm8.1}
\end{theorem}

   Define
$$
   F(t,x,\tau,\xi) = \frac{e^{-i\xi\cdot x - i \tau t}}{2\vert\xi\vert} 
	   \left[\frac{1 - e^{it(\tau + \vert \xi\vert)}}{\tau + \vert \xi\vert} - \frac{1 - e^{it(\tau - \vert \xi\vert)}}{\tau - \vert \xi\vert} \right].
$$
The next proposition gives the harmonizable representation of $\hat v$. This representation also appears in \cite[Section 6]{balan}.

\begin{prop} Set
\begin{align*}
   v(t,x) &= \int_{\IR}\int_{\IR^k} W(d\tau,d\xi)\, \vert \xi \vert^{(\beta - k)/2}\, F(t,x,\tau,\xi),
\end{align*}
and let $\tilde W_j(\varphi)$ be defined as in Proposition \ref{prop20}. Then $(v(t,x),\, (t,x) \in \IR_+ \times \IR^k)$ is a $\IcC$-valued solution of
$$
\left\lbrace
\begin{array}{rcl}
\frac{\partial^2}{\partial t^2} v_{j}(t,x) & = & \Delta v_{j}(t, x) + \dot {\tilde W}_{j} (t,x),\qquad j= 1, \dots, d,\\ [8pt]
v (0, x) & = & 0, \qquad \frac{\partial}{\partial t} v(0,x) = 0, \qquad x \in \IR.
\end{array}
\right.
$$
In particular, Re$(v)$ and $\hat v$ have the same law.
\label{prop7.2}
\end{prop}

\proof Let $S(s,y)$ be the fundamental solution of the wave equation. Since $\beta \in \, ]0,k\wedge2[$ or $k=1=\beta$, the stochastic integral
$$
   \int_{\IR}\int_{\IR^k} \tilde W_j(ds,dy) 1_{[0,t]}(s) S(t-s,x-y)
$$
is well-defined in all spatial dimensions $k\geq 1$ (see \cite[Example 6]{D99}), and
\begin{align*}
   & \int_{\IR}\int_{\IR^k} \tilde W_j(ds,dy) 1_{[0,t]}(s) S(t-s,x-y)\\
	 &\qquad = \int_{\IR}\int_{\IR^k} W_j(d\tau,d\xi)\, \vert \xi \vert^{(\beta - k)/2} \, \cF_{s,y}(1_{[0,t]}(\cdot) S(t-\cdot, x - \cdot))(\tau,\xi).
\end{align*}
Now for $s \in [0,t]$, according to \cite[Example 6]{D99},
$$
   \cF_y S(t-s,x-\cdot)(\xi) = e^{-i\xi\cdot x} \cF_y S(t-s,\cdot)(-\xi) = e^{-i\xi\cdot x} \frac{\sin((t-s)\vert \xi\vert)}{\vert \xi\vert},
$$
and $\cF_{s,y}(1_{[0,t]}(\cdot) S(t-\cdot, x - \cdot))(\tau,\xi)$ is equal to
\begin{align*}
	&\frac{e^{-i\xi\cdot x}}{\vert \xi\vert} \int_0^t e^{-i\tau s} \sin((t-s)\vert \xi\vert)\, ds 
	= \frac{e^{-i\xi\cdot x}}{\vert \xi\vert} \int_0^t e^{-i\tau (t-r)} \sin(r\vert \xi\vert) \, dr\\
	&\qquad = \frac{e^{-i\xi\cdot x-i\tau t}}{\vert \xi\vert} \int_0^t e^{i\tau r}\ \frac{e^{ir \vert \xi\vert} - e^{-ir\vert \xi\vert}}{2i}\, dr 
	= \frac{e^{-i\xi\cdot x-i\tau t}}{2\vert \xi\vert} \left[\frac{1-e^{it(\tau + \vert \xi\vert)}}{\tau + \vert \xi\vert} + \frac{e^{it(\tau - \vert \xi\vert)}-1}{\tau - \vert \xi\vert} \right].
\end{align*}
Therefore,
$$
   \int_{\IR}\int_{\IR} \tilde W_j(ds,dy) 1_{[0,t]}(s) S(t-s,x-y) = v_j(t,x),\qquad j=1,\dots,d.
$$
This proves the proposition.
\hfill $\Box$
\vskip 16pt

Let
$$
   \alpha = \frac{2 - \beta}{2},
$$
and set
$$
   v(A,t,x) = \int\!\!\!\int_{\max(\vert \tau\vert^\alpha, \vert \xi\vert^\alpha) \in A}  \vert \xi \vert^{(\beta - k)/2} F(t,x,\tau,\xi)\, W(d\tau,d\xi).
$$
Clearly, the random field $(v(A,t,x),\, A \in \ImB(\IR_+),\, (t,x) \in \IR_+ \times \IR^k)$ satisfies Assumption \ref{assump1}(a) (with the generic variable $x \in \IR^k$ replaced by $(t,x) \in \IR_+ \times \IR^k$). In the next lemma, we check Assumption \ref{assump1}(b).

\begin{lemma} Let $I \subset\, ]0,T]\times\IR$ be a compact box. Assume that \eqref{ebeta} holds. Then the random field $(v(A,t,x),\ A \in \ImB(\IR_+),\ (t,x) \in I)$ satisfies the conditions of Assumption \ref{assump1}, with exponents
$$
   \gamma_1 = \gamma_2 = \alpha^{-1} - 1 = \frac{\beta}{2 - \beta}=: \gamma.
$$
In particular, there is a universal constant $c_0$ and $a_0 \in \IR_+$ such that for all $a_0 \leq a \leq b$, $(t_0,x_0) \in I$, $(t,x) \in I$,
\begin{align}\nonumber
   &\Vert v([a,b[,t,x) - v(t,x) - v([a,b[,t_0,x_0) - v(t_0,x_0) \Vert_{L^2} \\
	&\qquad \leq c_0\left[a^\gamma \vert t - t_0 \vert + a^\gamma \sum_{j=1}^k \vert x_j - x_{0,j} \vert + b^{-1}\right]
	\label{eq8.1a}
\end{align}
and
\begin{align}
   \Vert v([0,a_0],t,x) - v([0,a_0],t_0,x_0) \Vert _{L^2} \leq c_0 \left[\vert t - t_0\vert + \sum_{j=1}^k  \vert x_j - y_j\vert \right].
		\label{eq8.1b}
\end{align}

\label{lem8.3}
\end{lemma}

\proof Assumption \ref{assump1}(a) is clearly satisfied, so we check Assumption \ref{assump1}(b). Let
$$
   D_1(a) = \{(\tau,\xi): \max(\vert\tau\vert^\alpha, \vert\xi\vert^\alpha) < a \}, \quad D_2(b) = \{(\tau,\xi): \max(\vert\tau\vert^\alpha, \vert\xi\vert^\alpha) > b \},
$$
and for $\ell=1,2$,
$$
   v_\ell(a,t,x) = \int\!\!\!\int_{D_\ell(a)} \vert \xi \vert^{(\beta - k)/2}\, F(t,x,\tau,\xi)\,  W(d\tau,d\xi).
$$
As in \eqref{vab}, 
\begin{align}\nonumber
 & v([a,b],t,x) - v(t,x) - v([a,b],t_0,x_0) + v(t_0,x_0) \\
 & \qquad = v_1(a,t_0,x_0) - v_1(a,t,x) + v_2(b,t_0,x_0) - v_2(b,t,x).
\label{eq8.1}
\end{align}
So for $\ell=1,2$, we let
$$
   f_\ell(a,t,x,t_0,x_0) = E\left[\vert v_\ell(a,t_0,x_0) - v_\ell(a,t,x)  \vert^2\right].
$$
Clearly,
$$
  f_1(a,t,x,t_0,x_0) = d \int\!\!\!\int_{D_1(a)} \vert F(t_0,x_0,\tau,\xi) - F(t,x,\tau,\xi) \vert^2\, \vert \xi \vert^{\beta - k}\, d\tau d\xi.
$$
Using Lemma \ref{lem8.4}(a) below, we see that
\begin{align*}
  f_1(a,t,x,t_0,x_0) &\leq c\left(\vert t - t_0\vert^2 + \sum_{j=1}^k \vert x_j - x_{0,j} \vert^2 \right)\\
	&\quad \times \int\!\!\!\int_{D_1(a)} \left[\frac{1}{1+\frac{1}{4} (\tau+\vert\xi\vert)^2} + \frac{1}{1+\frac{1}{4} (\tau-\vert\xi\vert)^2}\right]\, \vert \xi \vert^{\beta - k}\, d\tau d\xi .
\end{align*}
Change to polar coordinates $r = \vert \xi \vert$ to see that the double integral is equal to
$$
   C \int\!\!\!\int_{0 \leq \max(\vert\tau\vert,r) < a^{1/\alpha},\, r>0}  \left[\frac{1}{1+\frac{1}{4} (\tau+r)^2} + \frac{1}{1+\frac{1}{4} (\tau-r)^2}\right]\, r^{\beta - 1}\, d\tau dr.
$$
Use the change of variables $u=(\tau +r)/2$, $v=  (\tau -r)/2$ to see that the double integral is equal to
\begin{align*}
   &2 \int\!\!\!\int_{\max(\vert u+v\vert,\,  u-v) \leq 2 a^{1/\alpha},\, u-v>0 } \left[\frac{1}{1+u^2} + \frac{1}{1+v^2}\right] (u-v)^{\beta - 1} dudv \\
	&\qquad \leq \int\!\!\!\int_{\max(\vert u\vert,\, \vert v\vert ) \leq 2 a^{1/\alpha} } \left[\frac{1}{1+u^2} + \frac{1}{1+v^2}\right] \vert u-v\vert^{\beta - 1} dudv. \\
\end{align*}
By Lemma \ref{lem8.5}(a) below, this is $ \leq C a^{\beta/\alpha} = C a^{2\gamma}$. We conclude that
\begin{equation}\label{eq8.2}
   \Vert v_1(a,t_0,x_0) - v_1(a,t,x) \Vert_{L^2} \leq C\, a^{\gamma} \left[\vert t - t_0\vert + \sum_{j=1}^k \vert x_j - x_{0,j} \vert\right].
\end{equation}
This establishes in particular \eqref{eq8.1b}, for any $a_0 \geq 0$.

   We now turn to the second term:
\begin{align}\nonumber
   f_2(b,t,x,t_0,x_0) &= d \int\!\!\!\int_{D_2(b)} \vert F(t_0,x_0,\tau,\xi) - F(t,x,\tau,\xi) \vert^2\, \vert \xi \vert^{\beta - k}\, d\tau d\xi\\
	 &\leq 2 d \int\!\!\!\int_{D_2(b)} \left[(F(t_0,x_0,\tau,\xi))^2 + (F(t,x,\tau,\xi))^2\right]\, \vert \xi \vert^{\beta - k}\, d\tau d\xi.
\label{eq8.3}
\end{align}
By Lemma \ref{lem8.4}(b) below, the double integral is bounded above by
$$
  C_T \int\!\!\!\int_{D_2(b)} \left[\frac{1}{1+\frac{1}{4} \vert \tau + \vert\xi\vert \vert^2} + \frac{1}{1+\frac{1}{4} 
  \vert \tau - \vert\xi\vert \vert^2}\right] \frac{\vert \xi \vert^{\beta - k}}{1+ \vert\xi\vert^2} \, d\tau d\xi.
$$
Change again to polar coordinates $r = \vert\xi\vert$ to see that this is bounded by
\begin{equation}\label{eq8.4}
   \int\!\!\!\int_{\max(\vert \tau\vert ,r) > b^{1/\alpha},\, r>0} \left[\frac{1}{1+\frac{1}{4} (\tau+r)^2} + \frac{1}{1+\frac{1}{4} (\tau-r)^2}\right]\, \frac{r^{\beta - 1}}{1+ r^2}\, d\tau dr.
\end{equation}
By Lemma \ref{lem8.5}(b) below, this is $\leq c b^{-2}$.

   We conclude from \eqref{eq8.3} and the above estimate \eqref{eq8.4} that for large $b$,
\begin{equation}\label{eq8.5}
   f_2(b,t,x,t_0,x_0) \leq Cb^{-2}.
\end{equation}
Putting together \eqref{eq8.1}, \eqref{eq8.2} and \eqref{eq8.5}, we conclude that for $a_0$ large enough and $a_0 \leq a \leq b$, the conclusion of Lemma \ref{lem8.3} holds.
\hfill $\Box$
\vskip 16pt

   The following two lemmas were used in the proof of Lemma \ref{lem8.3}.
	
\begin{lemma} Fix $T > 0$. There is a constant $C_T$ such that for all $(t,x),\, (t_0,x_0) \in [0,T] \times \IR^k$, and all $(\tau,\xi) \in \IR \times \IR^k$, the following inequalities hold:

   (a)
\begin{align*}
  & \vert F(t_0,x_0,\tau,\xi) - F(t,x,\tau,\xi)\vert \\
	&\qquad \leq C_T \left(\vert t - t_0\vert +\sum_{j=1}^k\vert x_j - x_{0,j}\vert\right) 
	\left[\frac{1}{1+\frac{1}{2} \vert\tau+\vert\xi\vert\vert} + \frac{1}{1+\frac{1}{2} \vert\tau-\vert\xi\vert\vert}\right] .
\end{align*}

   (b)
$$
  \vert F(t,x,\tau,\xi)\vert \leq C_T\left[\frac{1}{1+\frac{1}{2} \vert \tau + \vert\xi\vert \vert} + \frac{1}{1+\frac{1}{2} 
  \vert \tau - \vert\xi\vert \vert}\right] \frac{1}{1+ \vert\xi\vert}.
$$
\label{lem8.4}
\end{lemma}

\proof (a) Notice that
\begin{align*}
   \frac{\partial F}{\partial x_j}(t,x,\tau,\xi) = -i\xi_j F(t,x,\tau,\xi) 
	= \frac{-i \xi_j}{2 \vert \xi\vert} e^{-i\xi\cdot x - i \tau t} \left[\frac{1 - e^{it(\tau + \vert \xi\vert)}}{\tau + \vert \xi\vert} - \frac{1 - e^{it(\tau - \vert \xi\vert)}}{\tau - \vert \xi\vert} \right].
\end{align*}
Observe that there is $c>0$ such that for all $u \in \IR$ and $t \in [0,T]$,
$$
   \left\vert \frac{1-e^{itu}}{u} \right\vert \leq \frac{c}{1+ \frac{1}{2} \vert u \vert},
$$
so
\begin{equation}\label{eq8.6}
   \left\vert \frac{\partial F}{\partial x_j}(t,x,\tau,\xi) \right\vert \leq \frac{c}{2} \left[\frac{1}{1+\frac{1}{2} \vert\tau+\vert\xi\vert\vert} + \frac{1}{1+\frac{1}{2} \vert\tau-\vert\xi\vert\vert}\right].
\end{equation}

   Similarly,
\begin{align*}
   \frac{\partial F}{\partial t}(t,x,\tau,\xi) & = -i\tau F(t,x,\tau,\xi) + \frac{e^{-i\xi\cdot x - i \tau t}}{2 \vert \xi\vert} \left[-ie^{it(\tau+\vert\xi\vert)} + i e^{it(\tau-\vert\xi\vert)} \right]\\
	&= \frac{-i}{2 \vert \xi\vert} e^{-i\xi\cdot x} \left[\frac{\tau e^{- i \tau t} -\tau e^{it\vert\xi\vert}}{\tau+\vert\xi\vert} - \frac{\tau e^{- i \tau t} - \tau e^{-it\vert\xi\vert}}{\tau-\vert\xi\vert} + e^{it\vert\xi\vert} - e^{-it\vert\xi\vert} \right]. 
\end{align*}
We notice that the term in brackets vanishes when $\vert\xi\vert =0$, and remains bounded when $\tau\pm \vert\xi\vert \to 0$, so $\frac{\partial F}{\partial t}$ is locally bounded. In fact, reducing to a common denominator, rearranging terms and simplifying, one finds that
$$
\frac{\partial F}{\partial t}(t,x,\tau,\xi) = \frac{i}{2} e^{-i\xi\cdot x - i \tau t} \left[\frac{1-e^{it(\tau-\vert\xi\vert)}}{\tau-\vert\xi\vert} + \frac{1-e^{it(\tau+\vert\xi\vert)}}{\tau+\vert\xi\vert} \right],
$$
therefore, as in \eqref{eq8.6},
\begin{equation}\label{eq8.7}
   \left\vert \frac{\partial F}{\partial t}(t,x,\tau,\xi) \right\vert \leq \frac{c}{2} \left[\frac{1}{1+\frac{1}{2} \vert\tau+\vert\xi\vert\vert} + \frac{1}{1+\frac{1}{2} \vert\tau-\vert\xi\vert\vert}\right].
\end{equation}
Using \eqref{eq8.6}, \eqref{eq8.7} and the Mean Value Theorem, we see that (a) holds.

   (b) Let $u=(\tau +\vert\xi\vert)/2$, $v=  (\tau -\vert\xi\vert)/2$ and notice that
\begin{equation}
   \vert F(t,x,\tau,\xi)\vert = \frac{\vert \varphi_t(2u) - \varphi_t(2v) \vert}{2 \vert u-v \vert},
\end{equation}
where
$$
   \varphi_t(u) = \frac{1 - e^{itu}}{u}, \qquad u\neq 0.
$$
Setting $\varphi_t(0) = -it$, then $\varphi_t \in C^1(\IR, \IcC)$, and
$$
   \varphi'_t(u) = \frac{-1 + e^{itu} - itu e^{itu}}{u^2} \qquad \mbox{if } u \neq 0,
$$
and $\varphi'_t(0) = t^2/2$. It follows that for all $(t,u) \in [0,T]\times \IR$,
\begin{equation}\label{eq8.8a}
  \max (\vert  \varphi_t(u)\vert, \vert  \varphi'_t(u)\vert) \leq \frac{C_T}{1+ \vert u \vert}.
\end{equation}
In particular, we claim that for all $(t,u) \in [0,T]\times \IR$ with $\vert u - v \vert \leq 1/2$,
\begin{equation}\label{eq8.9}
   \frac{\vert \varphi_t(2u) - \varphi_t(2v) \vert}{\vert u-v \vert} \leq C_T\left[\frac{1}{1+ \vert u \vert} + \frac{1}{1+ \vert v \vert} \right].
\end{equation}
Indeed, by the Mean Value Theorem,
$$
   \vert \varphi_t(2u) - \varphi_t(2v) \vert \leq 2 \vert u-v \vert \, \vert \varphi'_t(\xi) \vert ,
$$
for some $\xi$ between $u$ and $v$. If both $u$ and $v$ have the same sign, say if $0<u<v$, then by \eqref{eq8.8a},
$$
   \vert \varphi'_t(\xi) \vert \leq \frac{C}{1+ \vert \xi \vert} \leq \frac{C}{1+ \vert u \vert} \leq C\left[ \frac{1}{1+ \vert u \vert} + \frac{1}{1+ \vert v \vert}\right].
$$
The case where $u$ and $v$ are both negative is handled similarly. Finally, if $u<0<v$, then since $\vert u - v \vert \leq 1/2$, we have $\vert u \vert \leq 1/2$ and $\vert v \vert \leq 1/2$, so
$$
   \vert \varphi'_t(\xi) \vert \leq \frac{C}{1+ \vert \xi \vert} \leq C = C \left[\frac{3/4}{1+ \frac{1}{2}} + \frac{3/4}{1+ \frac{1}{2}}\right] \leq \tilde C \left[\frac{1}{1+ \vert u \vert} + \frac{1}{1+ \vert v \vert} \right].
$$
This proves \eqref{eq8.9}.

   We now claim that there is a constant $C_T < \infty$ such that for all $(t,u,v) \in \IR_+ \times \IR^2$,
\begin{equation}\label{eq7.12}
    \frac{\vert \varphi_t(2u) - \varphi_t(2v) \vert}{2 \vert u-v \vert} \leq C_T\left[\frac{1}{1+ \vert u \vert} + \frac{1}{1+ \vert v \vert} \right] \frac{1}{1+ \vert u-v\vert}.
\end{equation}
Indeed, assume first that $\vert u-v\vert \leq 1/2 $. Then by \eqref{eq8.9}, the left-hand side is
$$
   \leq C_T \left[\frac{1}{1+ \vert u \vert} + \frac{1}{1+ \vert v \vert} \right] \leq C_T \left[\frac{1}{1+ \vert u \vert} + \frac{1}{1+ \vert v \vert} \right] \frac{3/2}{1+ \vert u-v\vert}.
$$
Now assume that $\vert u-v\vert \geq 1/2 $. Then the left-hand side of \eqref{eq7.12} is
$$
   \leq \frac{3/2}{1+ \vert u-v\vert} (\vert \varphi_t(2u)\vert + \vert\varphi_t(2v) \vert) \leq \frac{\tilde C_T}{1+ \vert u-v\vert} \left[ \frac{1}{1+ \vert u \vert} + \frac{1}{1+ \vert v \vert}\right],
$$
where we have used \eqref{eq8.8a}. This completes the proof of (b).
\hfill $\Box$
\vskip 12pt

\begin{lemma} (a) For $\beta \in \, ]0,2[$, 
$$
   \int\!\!\!\int_{\max(\vert u\vert,\, \vert v\vert ) \leq 2 a^{1/\alpha} } \left[\frac{1}{1+u^2} + \frac{1}{1+v^2}\right] \vert u-v\vert^{\beta - 1} dudv \leq C a^{\beta /\alpha}.
$$

   (b) If $\beta \geq 1$, then for large $b$,
$$
   \int\!\!\!\int_{\max(\vert \tau\vert ,r) > b^{1/\alpha},\, r>0} \left[\frac{1}{1+\frac{1}{4} (\tau+r)^2} + \frac{1}{1+\frac{1}{4} (\tau-r)^2}\right]\, \frac{r^{\beta - 1}}{1+ r^2}\, d\tau dr \leq C b^{-2}.
$$
\label{lem8.5}
\end{lemma}

\proof (a) It suffices to consider the two integrals
\begin{align*}
   A_1 &= \int\!\!\!\int_{\max( u,\,  v ) \leq 2 a^{1/\alpha} ,\, u>0,\, v>0} \left[\frac{1}{1+u^2} + \frac{1}{1+v^2}\right] \vert u-v\vert^{\beta - 1} dudv, \\
	 A_2 &= \int\!\!\!\int_{\max( u,\,  v ) \leq 2 a^{1/\alpha} ,\, u>0,\, v>0} \left[\frac{1}{1+u^2} + \frac{1}{1+v^2}\right] \vert u+v\vert^{\beta - 1} dudv. 
\end{align*}
By symmetry, $A_1 = 2A_{1,1}$, where
\begin{align*}
  A_{1,1} &= \int_0^{2a^{1/\alpha}} du \int_0^u dv \left[\frac{1}{1+u^2} + \frac{1}{1+v^2}\right] ( u-v)^{\beta - 1} dudv \\
	   &\leq 2 \int_0^{2a^{1/\alpha}} du \int_0^u dv\, \frac{1}{1+v^2} (u-v)^{\beta - 1}.
\end{align*}
By Fubini's theorem, this is equal to
\begin{align*}
  & 2 \int_0^{2a^{1/\alpha}} \frac{dv}{1+v^2} \int_v^{2a^{1/\alpha}} du\, (u-v)^{\beta-1} \\
	&\qquad = 2 \int_0^{2a^{1/\alpha}} \frac{dv}{1+v^2}\, (2a^{1/\alpha} -v)^\beta 
\leq C \int_0^{2a^{1/\alpha}} \frac{dv}{1+v^2}\, (2a^{1/\alpha})^\beta 
\leq C a^{\beta/\alpha} \int_0^\infty  \frac{dv}{1+v^2}.
\end{align*}

   Turning to $A_2$, by symmetry,
\begin{align*}
   A_2 & = 2 \int_0^{2 a^{1/\alpha}} dv \int_0^v du \left[\frac{1}{1+u^2} + \frac{1}{1+v^2}\right] \vert u+v\vert^{\beta - 1} \\
	& \leq C \int_0^{2 a^{1/\alpha}} dv \int_0^v du\, \frac{1}{1+u^2}\, (u+v)^{\beta - 1}.
\end{align*}
By Fubini's Theorem, this is equal to
\begin{align*}
   C \int_0^{2 a^{1/\alpha}}  \frac{du}{1+u^2} \int_u^{2 a^{1/\alpha}} dv\, (u+v)^{\beta - 1} &=  C \int_0^{2 a^{1/\alpha}}  \frac{du}{1+u^2} \left[(2 a^{1/\alpha} +u)^\beta - (2u)^\beta \right]\\
	&\leq C \int_0^{2 a^{1/\alpha}}  \frac{du}{1+u^2}\, (3 a^{1/\alpha})^\beta \leq  \tilde C a^{\beta/\alpha}.
\end{align*}
This proves (a).

   (b) We need to integrate over two regions:
$$
   r > b^{1/\alpha},\ \vert\tau\vert < r, \qquad\mbox{and}\qquad \vert\tau\vert > b^{1/\alpha},\ 0 < r < \vert\tau\vert.
$$
Concerning the first region, we have to consider
$$
   \int_{b^{1/\alpha}}^\infty dr \int_{-r}^r d\tau\, \frac{1}{1+\frac{1}{4} (\tau\pm r)^2} \frac{r^{\beta - 1}}{1+ r^2},
$$
and, by symmetry, it suffices to consider
$$
   \int_{b^{1/\alpha}}^\infty dr\, \frac{r^{\beta - 1}}{1+ r^2} \int_0^r d\tau \, \frac{1}{1+\frac{1}{4} (\tau\pm r)^2} \leq \frac{\pi}{2} \int_{b^{1/\alpha}}^\infty dr\, r^{\beta - 3} = b^{(\beta - 2)/\alpha} = b^{-2}.
$$
For the second region, we consider
\begin{align*}
    \int_{b^{1/\alpha}}^\infty d\tau \int_{0}^\tau dr\, \frac{1}{1+\frac{1}{4} (\tau\pm r)^2} \frac{r^{\beta - 1}}{1+ r^2} &= \int_0^{b^{1/\alpha}} dr\, \frac{r^{\beta - 1}}{1+ r^2} \int_{b^{1/\alpha}}^\infty d\tau \, \frac{1}{1+\frac{1}{4} (\tau\pm r)^2} \\
		&\quad + \int_{b^{1/\alpha}}^\infty dr\, \frac{r^{\beta - 1}}{1+ r^2} \int_r^\infty d\tau \, \frac{1}{1+\frac{1}{4} (\tau\pm r)^2} .
\end{align*}

   The second integral is
$$
   \leq \int_{b^{1/\alpha}}^\infty dr\, r^{\beta - 3} \frac{\pi}{2} \leq C b^{(\beta - 2)/\alpha} = C b^{-2}.
$$
Concerning the first integral, in the case of a ``$+$" sign, it is
$$
 \leq \int_0^{b^{1/\alpha}} dr\, \frac{r^{\beta - 1}}{1+ r^2} \int_{b^{1/\alpha}}^\infty d\tau \, \frac{1}{1+\frac{1}{4} \tau^2} = \int_0^{b^{1/\alpha}} dr\, \frac{r^{\beta - 1}}{1+ r^2} \left[\frac{\pi}{2} - \arctan\left(\frac{b^{1/\alpha}}{2}\right)\right].
$$
Using the property
$$
   \lim_{x \to \infty} x \left[\frac{\pi}{2} - \arctan(x)\right] = 1,
$$
we see that for all $b \geq 1$, this is
$$
   \leq \tilde c\, b^{-1/\alpha} \int_0^{\infty} dr\, \frac{r^{\beta - 1}}{1+ r^2} \leq C b^{-2},
$$
since $\frac{1}{\alpha} = \frac{2}{2 - \beta}  \geq 2$ because $\beta \geq 1$.

   In the case of a ``$-$" sign, we write the first integral as $I_1(b) + I_2(b)$, where
\begin{align*}
   I_1(b) &= \int_0^{b^{1/\alpha}/2} dr\, \frac{r^{\beta - 1}}{1+ r^2} \int_{b^{1/\alpha}}^\infty d\tau \, \frac{1}{1+\frac{1}{4} (\tau- r)^2},\\
	 I_2(b) &= \int_{b^{1/\alpha}/2}^{b^{1/\alpha}} dr\, \frac{r^{\beta - 1}}{1+ r^2} \int_{b^{1/\alpha}}^\infty d\tau \, \frac{1}{1+\frac{1}{4} (\tau- r)^2}.
\end{align*}
Then
\begin{align*}
   I_1(b) &= \int_0^{b^{1/\alpha}/2} dr\, \frac{r^{\beta - 1}}{1+ r^2} \int_{b^{1/\alpha}-r}^\infty du\, \frac{1}{1+\frac{1}{4} u^2} \\
	 & \leq \int_0^{b^{1/\alpha}/2} dr\, \frac{r^{\beta - 1}}{1+ r^2} \int_{b^{1/\alpha}/2}^\infty du\, \frac{1}{1+\frac{1}{4} u^2} \leq C b^{-1/\alpha} \int_0^{b^{1/\alpha}/2} dr\, \frac{r^{\beta - 1}}{1+ r^2} \leq c b^{-2},
\end{align*}
since $\beta \geq 1$, and
\begin{align*}
I_2(b) & \leq \int_{b^{1/\alpha}/2}^{b^{1/\alpha}} dr\, \frac{r^{\beta - 1}}{1+ r^2} \int_{-\infty}^\infty d\tau \, \frac{1}{1+\frac{1}{4} (\tau- r)^2} \leq C \int_{b^{1/\alpha}/2}^{b^{1/\alpha}} dr\, r^{\beta - 3} \\
& \leq C b^{(\beta-2)/\alpha} = C b^{-2}.
\end{align*}
This completes the proof of (b).
\hfill $\Box$
\vskip 16pt

   We now turn to Assumption \ref{assump2}. In the context of this section, in agreement with Lemma \ref{lem8.3} and Assumption \ref{assump1}(b),
$$
	 \Delta((t,x),(s,y)) = \vert t-s \vert^{\frac{2-\beta}{2}} + \vert x-y \vert^{\frac{2-\beta}{2}}.
$$
It is well-known (see \cite[Proposition 1.4]{DSb}) that for any compact box $I \subset \,]0,\infty[ \times \IR^k$, there is $c>0$ such that
$$
    \Vert \hat v(t,x) - \hat v(s,y) \Vert_{L^2}^2 \geq c\, \Delta((t,x),(s,y)).
$$
Further, using the change of variables $r=t-s$, $\eta = (t-s) \xi$, we see that
\begin{align*}
    \Vert \hat v(t,x) \Vert_{L^2}^2 &= d \int_0^t ds \int_{\IR^k} \frac{d\xi}{\vert \xi \vert^{k-\beta}} \, \frac{\sin^2((t-s) \vert \xi \vert)}{\vert \xi \vert^2} 
		= d \int_0^t dr\, r^{2-\beta} \int_{\IR^k} \frac{d\eta}{\vert \eta \vert^{k+2-\beta}} \sin^2(\vert \eta \vert) 
		\\ &
		= c t^{3-\beta},
\end{align*}
so Assumption \ref{assump2}(a) is satisfied for the box $I$. In the next lemma, we check Assumption \ref{assump2}(b).

\begin{lemma} Let $I \subset \, ]0,\infty[ \times \IR^k$ be a compact box. Fix $(t,x) \in I$. Let $t' = t - 2(2\rho)^{\alpha^{-1}}$ and $x'=x$ (where $\rho$ is small enough so that $t' >0$). Assume that $k=1 = \beta$ or $1 < \beta < k \wedge 2$. There is a number $C_1$ (depending on $\rho$, $\beta$, $k$ and $d$) such that for all $(s_1,y_1), (s_2,y_2) \in B_\rho'(t,x)$ (the open $\Delta$-ball of radius $2\rho$ centered at $(t,x)$) and $j\in \{1,\dots,d\}$,
\begin{equation}\label{e7.13}
   \left\vert E[(\hat v_j(s_1,y_1) - \hat v_j(s_2,y_2)) \hat v_j(t',x')]\right\vert  \leq C_1 \left(\vert s_1 - s_2 \vert^{\delta} +  \vert y_1 - y_2 \vert^{\delta}\right),
\end{equation}
where $\delta = 2 - \beta$.
\label{lem7.6}
\end{lemma}

\begin{remark} The conclusion of this lemma would not be possible for $\beta \in \,]0,1[$, since it would mean that $(s,y) \mapsto E[(\hat v_j(s,y) \hat v_j(t',x')]$ would be H\"older-continuous with exponent $2 - \beta >1$.
\end{remark}

\proof Consider first the case $k=1=\beta$ (space-time white noise in spatial dimension $k=1$). Then
\begin{equation}\label{e7.14}
   E[(\hat v_j(s,y) \hat v_j(t',x')] = \int_0^{t'} dr \int_\IR dz \, S(s-r,y-z) S(t'-r,x'-z),
\end{equation}
where $S(r,z) = \frac{1}{2} 1_{\{\vert z \vert < r \}}$ is the fundamental solution of the wave equation. For $(s,y) \in B_\rho'(t,x)$,
$$
   \vert s - t \vert \leq (2\rho)^{\alpha^{-1}}, \qquad  \vert y - x \vert \leq (2\rho)^{\alpha^{-1}},
$$
and since $t' = t - 2(2\rho)^{\alpha^{-1}}$ and $x'=x$, one checks immediately that if, in addition,
$
   \vert x' - z \vert < t' - r, 
$
then
$$
   \vert y - z\vert \leq \vert y - x\vert + \vert x - z\vert < (2\rho)^{\alpha^{-1}} + t'-r = t - (2\rho)^{\alpha^{-1}} - r \leq s-r,
$$
so the right-hand side of \eqref{e7.14} is equal to
$$
    \int_0^{t'} dr \int_\IR dz \, \frac{1}{4} \, 1_{\{\vert x'-z\vert < t'-r\}},
$$
and therefore $(s,y) \mapsto E[(\hat v_j(s,y) \hat v_j(t',x')]$ is constant over $B_\rho'(t,x)$ and \eqref{e7.13} is trivially satisfied.

   Now consider the case where $1 < \beta < k \wedge 2$. Then for $s \geq t'$,
\begin{align*}
  & E[(\hat v_j(s,y) \hat v_j(t',x')] \\
	&\qquad = \int_0^{t'} dr \int_{\IR^k} \frac{d\xi}{\vert \xi \vert^{k-\beta}} \cF S(s-r,y - \cdot)(\xi) \overline{\cF S(t'-r,x' - \cdot)(\xi)} \\
	 &\qquad = \int_0^{t'} dr \int_{\IR^k} d\xi\, \vert \xi \vert^{\beta - 2 - k} e^{-i\xi\cdot(y-x')} \sin((s-r) \vert \xi \vert) \sin((t'-r)\vert \xi \vert) \\
	&\qquad = \int_0^{t'} dr \int_{\IR^k} d\xi\, \vert \xi \vert^{\beta - 2 - k} e^{-i\xi\cdot(y-x')} \sin(r \vert \xi \vert) \sin((h+r)\vert \xi \vert),
\end{align*}
where we have set $h=s-t'$. We now permute the two integrals and calculate the $dr$-integral explicitly. As in the proof of Lemma A.12 in \cite{DSm}, this gives
\begin{equation}\label{e7.15}
   E[\hat v_j(s,y) \hat v_j(t',x')] = (t')^{3-\beta} \int_{\IR^k} d\eta \, \vert \eta \vert^{\beta - 2 - k} e^{-i\eta\cdot u} g_0(\lambda, \vert \eta \vert),
\end{equation}
where $u=(y-x')/t'$, $\lambda = (s-t')/t'$, and
$$
   g_0(\lambda,r) = \cos(\lambda r) - \frac{\sin(r)}{r}\, \cos((\lambda + 1) r).
$$
\vskip 12pt

\noindent{\em Case 1} (time increments): $s_1 \neq s_2$, $y_1 = y_2 = y$. Set $\lambda_1 = (s_1-t')/t'$, $\lambda_2 = (s_2-t')/t'$, and $u = (y-x')/t'$. Then by \eqref{e7.15}, 
\begin{align}\nonumber
   & E[(\hat v_j(s_1,y) - \hat v_j(s_2,y)) \hat v_j(t',x')] \\ 
	 &\qquad =  (t')^{3-\beta} \int_{\IR^k} d\eta \, \vert \eta \vert^{\beta - 2 - k} e^{-i\eta\cdot u} (g_0(\lambda_1, \vert \eta \vert) - g_0(\lambda_2, \vert \eta \vert)).
\label{e7.16}
\end{align}
Because $\cos(\cdot)$ is Lipschitz, we see that
$$
   \vert g_0(\lambda_1, r) - g_0(\lambda_2, r) \vert \leq 4( (\vert \lambda_1- \lambda_2\vert r) \wedge 1),
$$
therefore, the right-hand side of \eqref{e7.16} is bounded above by $4 (t')^{3-\beta} (I_1 + I_2)$, where
\begin{align*}
   I_1  = \vert \lambda_1- \lambda_2\vert \int_0^{\vert \lambda_1- \lambda_2\vert^{-1}}  r^{\beta-2}\, dr, \qquad
	 I_2  = \int_{\vert \lambda_1- \lambda_2\vert^{-1}}^{\infty}   r^{\beta-3}\, dr.
\end{align*}
Clearly, since $\beta > 1$,
$$
   I_1 = \vert \lambda_1- \lambda_2\vert \, \frac{\vert \lambda_1- \lambda_2\vert^{1-\beta}}{\beta - 1} = c\, \vert s_1-s_2 \vert^{2-\beta},
$$
and
$$
   I_2  = \frac{\vert \lambda_1- \lambda_2\vert^{2-\beta}}{2-\beta} = c\, \vert s_1-s_2 \vert^{2-\beta}.
$$
We conclude that
\begin{equation}\label{e7.17}
   E[(\hat v_j(s_1,y) - \hat v_j(s_2,y)) \hat v_j(t',x')] \leq c\, \vert s_1-s_2 \vert^{2-\beta}.
\end{equation}
\vskip 12pt

\noindent{\em Case 2} (spatial increments): $s_1 = s_2 = s$, $y_1 \neq y_2$. Set $\lambda = (s-t')/t'$, $u_1 = (y_1 - x')/t'$, $u_2 = (y_2 - x')/t'$. By \eqref{e7.15}, 
\begin{align}
	E[(\hat v_j(s,y_1) - \hat v_j(s,y_2)) \hat v_j(t',x')] 
	=  (t')^{3-\beta} \int_{\IR^k} d\eta \, \vert \eta \vert^{\beta - 2 - k} (e^{-i\eta\cdot u_1} - e^{-i\eta\cdot u_2}) g_0(\lambda, \vert \eta \vert).
\label{e7.18}
\end{align}
Notice that
$$
   \vert e^{-i\eta\cdot u_1} - e^{-i\eta\cdot u_2} \vert \leq 2 ((\vert u_1 - u_2 \vert\, \vert \eta \vert) \wedge 1)
$$
and $g_0(\lambda,r) \leq 2$, so the right-hand side of \eqref{e7.18} is bounded above by $4(t')^{3-\beta}(J_1 + J_2)$, where
$$
   J_1  = \vert u_1- u_2\vert \int_0^{\vert u_1- u_2\vert^{-1}}  r^{\beta-2}\, dr, \qquad
	 J_2  = \int_{\vert u_1- u_2\vert^{-1}}^{\infty}   r^{\beta-3}\, dr.
$$
Clearly, the same calculations as for $I_1$ and $I_2$ show that
$$
   J_1 + J_2 \leq \tilde c \, \vert u_1- u_2\vert^{2-\beta} = c \, \vert y_1- y_2\vert^{2-\beta}.
$$
We conclude that
\begin{equation}\label{e7.19}
   E[(\hat v_j(s,y_1) - \hat v_j(s,y_2)) \hat v_j(t',x')] \leq c \, \vert y_1- y_2\vert^{2-\beta}.
\end{equation}

   Putting together \eqref{e7.17} and \eqref{e7.19} establishes \eqref{e7.13}. This proves Lemma \ref{lem7.6}.
\hfill $\Box$
\vskip 16pt

\noindent{\em Proof of Theorem \ref{thm8.1}.} By Lemma \ref{lem8.3} and the sentences that precede this lemma, for any compact box $I \subset \,]0,\infty[\times \IR^k$, Assumption \ref{assump1} is satisfied for Re$(v)$, with exponents $\gamma_1 = \frac{\beta}{2-\beta} = \gamma_j$, $j=2,\dots,k+1$, so that $\alpha_1 = \frac{2-\beta}{2} = \alpha_j$, $j=2,\dots,k+1$. By Lemma \ref{lem7.6} and the comments that precede this lemma, Assumption \ref{assump2} is satisfied by $\hat v$ (with $\delta_j = 2-\beta > \alpha_j$), hence by Re$(v)$ by Proposition \ref{prop7.2}. Since $Q = \alpha_1^{-1} + k \alpha_2^{-1} = (2+2k)/(2-\beta) = d$, it follows from Theorem \ref{thm1} that for all $z \in \IR^Q$,
$$
   P\{\exists (t,x) \in I: \hat v(t,x) = z \} = P\{\exists (t,x) \in I: \mbox{Re}(v(t,x)) = z \} = 0.
$$
Since this holds for all compact boxes $I \subset \,]0,\infty[\times \IR^k$, Theorem \ref{thm8.1} is proved.
\hfill $\Box$
\vskip 16pt

\noindent{\sc Acknowledgment.} The research reported in this paper was initiated at the Centre Interfacultaire Bernoulli, Ecole Polytechnique F\'ed\'erale de Lausanne, Switzerland, during the semester program ``Stochastic Analysis and Applications" in Spring 2012. We thank this institution for its hospitality and support. 


\end{document}